\documentclass[letter,11pt]{amsart}
\usepackage{amssymb}
\newtheorem{thm}{Theorem}[section]
\newtheorem{lem}[thm]{\textbf Lemma}
\newtheorem{cor}[thm]{Corollary}

\theoremstyle{definition}

\newtheorem{defn}[thm]{Definition}

\theoremstyle{remark}

\newcommand{\les}{\lesssim}

\newcommand{\beeq}{\begin{equation}}
\newcommand{\eneq}{\end{equation}}

\newcommand{\supp}{\mbox{\rm supp}}

\newcommand{\eps}{{\varepsilon}}
\newcommand{\R}{{\mathbb R}}

\newcommand{\calL}{{\mathcal L}}
\newcommand{\calS}{{\mathcal S}}

\newcommand{\calH}{{\mathcal H}}

\newcommand{\sign}{\mbox{sign}}

\newcommand{\less}{\lesssim}
\newcommand{\gtr}{\gtrsim}
\newcommand\nind{\noindent}
\def\vol{{\rm vol}}
\def\la{\langle}
\def\ra{\rangle}
\begin{document}
\numberwithin{equation}{section}

\title[Decay for the Schr\"odinger evolution]
{Decay estimates for the Schr\"odinger evolution on asymptotically
conic surfaces of revolution I}
\author{Wilhelm Schlag, Avy Soffer, Wolfgang Staubach}
\thanks{The first and second authors were partly supported by  the National
Science Foundation.}

\address{first and third authors: University of Chicago, Department of Mathematics,
5734 South University Avenue, Chicago, IL 60637, U.S.A.}
\email{schlag@math.uchicago.edu, wolf@math.uchicago.edu}

\address{second author: Rutgers University, Department of Mathematics, 110 Freylinghuysen Road, Piscataway, NJ 08854, U.S.A.}

\email{soffer@math.rutgers.edu}

\maketitle
\section{Introduction}

\noindent It is well-known that the free Schr\"odinger evolution
satisfies the dispersive bound
\begin{equation}
\label{eq:disp} \| e^{it\Delta} f\|_\infty \less |t|^{-\frac{d}{2}}
\|f\|_1
\end{equation}
where $\Delta$ denotes the Laplacean in $\R^d$. Another instance of
such decay bounds are the global Strichartz estimates

\beeq\label{eq:strich} \| e^{it\Delta}
f\|_{L^{2+\frac{4}{d}}(\R^{d+1})} \less \|f\|_{L^2(\R^d)} \eneq

\nind and mixed-norm versions thereof. In this paper we establish a
decay estimate (valid for all $t$), similar to \eqref{eq:disp}, for
the Schr\"odinger evolution on a non-compact 2-dimensional manifold
with a trapped geodesic. As we shall explain below, the case of the
manifold considered in this investigation, as well as the method of
proving the decay estimate, are different from the studies in the
existing literature concerning Schr\"odinger evolution on manifolds.\\

There has been much activity lately around establishing dispersive
and Strichartz estimates for more general operators, namely for
Schr\"odinger operators of the form $H=-\Delta+V$ with a decaying
potential $V$ or even more general perturbations. The seminal paper
here is Jorne\'e-Soffer-Sogge\cite{JSS}, and we refer the reader to
the survey~\cite{Sch2} for more recent references in this area.

Around the same time as \cite{JSS}, Bourgain\cite{Bou} found
Strichartz estimates on the torus. This is remarkable, as compact
manifolds do not exhibit dispersion as in~\eqref{eq:disp} which was
always considered a key ingredient for~\eqref{eq:strich}. The theme
of Strichartz estimates on manifolds (both local and global in time)
was then developed further in several important papers, see
Smith-Sogge\cite{SS}, Staffilani-Tataru\cite{ST},
Burq-Gerard-Tzvetkov\cite{BGT1}, \cite{BGT2},
Hassel-Tao-Wunsch\cite{HTW1}, \cite{HTW2}, Robbiano-Zuily\cite{RZ},
and Tataru\cite{T}. Gerard\cite{G} reviews some of the recent work
in this field.

A recurring theme in this area is the importance of closed geodesics
for Strichartz estimates. In fact, it is well-known that the
presence of closed geodesics  necessarily leads to a loss of
derivatives in the Strichartz bounds. The intuition here is that
initial data that are highly localized around a closed geodesic and
possess high momentum traveling around this geodesic will lead to
so-called meta-stable states in the Schr\"odinger evolution. These
are states that remain "coherent" for a long time, which amounts to
absence of dispersion during that time, see for example~\cite{G} (in
the classical approximation, dispersive estimates are governed by
the Newtonian scattering trajectories - classically speaking, closed
geodesics are  non-scattering states).

For this reason, many authors have imposed explicit non-trapping
conditions, see \cite{SS}, \cite{HTW1}, \cite{HTW2}, \cite{RodTao}.
The relevance of this condition lies with the construction of a
parametrix, which always involves solving for suitable
bi-characteristics. On manifolds these bi-characteristics are
governed by the geodesics flow in the co-tangent bundle - hence the
relevance of closed geodesics.

There is a large body of work on the so-called Kato smoothing
estimates where this non-trapping condition also features
prominently, see for example Craig-Kappeler-Strauss\cite{CKS},
Doi\cite{Doi}, and Rodnianski-Tao\cite{RodTao}.

Thus, the literature on the trapping case is very limited. In
addition, we are not aware of a reference that
studies~\eqref{eq:disp} rather than~\eqref{eq:strich} on manifolds,
which then necessarily need to be noncompact.\\

\noindent In this paper, we consider surfaces of revolution

\begin{equation*}
\calS = \{(x,r(x)\cos\theta,r(x)\sin\theta ):-\infty< x<\infty
,0\leq\theta\leq 2 \pi\}
\end{equation*}
with the metric $ds^2=r^2(x)d\theta^2+(1+r'(x)^2)dx^2$. Examples of
such surfaces abound, and there is no point to developing a theory
for all of them simultaneously (consider the cases: $r(x)$ constant,
$r(x)$  periodic, and $r(x)$ rapidly growing). Rather, we single out
a class of surfaces of revolution which have rather explicit
behavior at both ends. Our basic examples are $r(x)=\langle
x\rangle^{\alpha}$ with $0\le\alpha<\infty$. In fact, in order not
to obscure our ideas by technical details, we will present the main
result of this paper for the case $\alpha=1$. We remark, however,
that the method of this paper equally well applies to other values
of~$\alpha$.

We now define the class of asymptotically conic manifolds we shall
mainly work with.

\begin{defn}
\label{def:cones} We assume that $\inf_xr(x)>0$, and asymptotically
that

 \beeq\label{0} r(x)=|x|h(x)\qquad \mathrm{for}\; |x|\geq 1
\eneq

\noindent where $h(x)=1+O(x^{-2})$ and also $h^{(k)}(x)=O(x^{-2-k})$
for all $k\geq 1$.  Examples of such $r(x)$ are
$r(x)=\sqrt{1+|x|^2}=:\langle x\rangle$ and variants thereof.  Note
that our surface $\calS$ is asymptotic to cones at both ends. For
convenience, we shall also make the symmetry assumption
$r(x)=r(-x)$. We shall refer to such a surface of revolution as {\em
asymptotically conic and symmetric}.
\end{defn}

\nind The main result of this paper is the following theorem,
where $\Delta_{\calS}$ denotes the Laplace-Beltrami operator
on~$\calS$.

\begin{thm}\label{thm1}
Let $\calS$ be a surface of revolution in $\R^3$ which is
asymptotically conic and symmetric as in Definition~\ref{def:cones}.
Then for all $t$ \beeq\label{10} \Vert
e^{it\Delta_{\calS}}f\Vert_{L^{\infty}(\calS )}\les |t|^{-1}\Vert
f\Vert_{L^1(\calS )} \eneq \noindent provided $f$ does not depend on
the angular variable $\theta$.
\end{thm}

The symmetry assumption can be easily removed, but we include it to
simplify the exposition. Section~\ref{sec:proof} will be devoted to
the proof of this theorem. We remark that for surfaces which are
asymptotic to $\la x\ra^\alpha$. It is clear that in case $\alpha
=0$ the surface of revolution $\calS$ is just the cylinder
$S^{1}\times\mathbb{R}$, with the metric $ds^2=d\theta^2+dx^2$, for
which the dispersive estimate with radial data is the same as the
one for the free 1-dimensional Schr\"odinger operator, namely one
has $\Vert
e^{it\Delta_{(S^{1}\times\mathbb{R})}}f\Vert_{L^{\infty}(S^{1}\times\mathbb{R}
)}\les |t|^\frac{-1}{2}\Vert f\Vert_{L^{1}(\mathbb{R})}$. For
$0<\alpha<1$, the methods of this paper yield the decay rate
$t^{-\frac12(1+\alpha)}$, whereas for $1\le\alpha<\infty$ it is
$t^{-1}$. The intuition concerning these decay rates is as follows:
Let $B(p,t)$ denote the volume of a geodesic ball of radius $t$
centered at the point $p\in\calS$. Then for fixed $p$ one has
\[ \vol(B(p,t)) \sim t^{1+\alpha} \text{\ \ or\ \ }t^2 \]
depending on whether $\alpha<1$ or $\alpha\ge1$. In view of the unitarity of
the Schr\"odinger flow, we see that the decay rate should be given by $\vol(B(p,t))^{-\frac12}$
and this is indeed the case.

In a sequel to this paper we will discuss the case of non-radial
data $f=f(\theta, r)$. In fact, the methods of this paper allow us
to prove the following result.

\begin{thm}\label{thm2}
Let $\calS$ be a surface of revolution in $\R^3$ which is
asymptotically conic and symmetric as in Definition~\ref{def:cones}.
Then for each integer $n$ there exists a constant $C(n)$ so that for
all~$t$ \beeq\nonumber \Vert
e^{it\Delta_{\calS}}(e^{in\theta}f)\Vert_{L^{\infty}(\calS )}\le
C(n)\, |t|^{-1}\Vert f\Vert_{L^1(\calS )} \eneq \noindent provided
$f$ does not depend on the angular variable $\theta$.
\end{thm}

Summing in $n$ we of course obtain a global $L^1(\calS)\to
L^\infty(\calS)$ decay estimate with a loss of derivatives in
$\theta$. The details, as well as the dependance of the constant
$C(n)$ on $n$ will be discussed in the sequel to this paper.
Obviously, the behavior of $C(n)$ for large~$n$ is very important as
it governs how many derivatives (in~$\theta$) we will lose on the
right-hand side. Note that the loss of derivatives in $\theta$ is in
agreement with the aforementioned intuition that meta-stable states
can form from data with high momentum around a closed geodesic. We
do not lose any derivatives in Theorem~\ref{thm1} since radial data
cannot get trapped in the classical picture (the scattering
trajectories are generators of our surface of revolution).  Finally,
we remark that Theorem~\ref{thm1} leads to a Strichartz estimate for
radial data using standard techniques.

We now briefly describe the main ideas behind Theorems~\ref{thm1}
and~\ref{thm2}. First, using arc-length coordinates $\xi$ on $\calS$
and after multiplying by the weight $r^{\frac12}(\xi)$, we reduce
matters to the Schr\"odinger operator (with $n$ as in
Theorem~\ref{thm2})
\[ \calH_n:=-\partial_\xi^2 + V(\xi) + \frac{n^2}{r^2(\xi)} \]
on the line.  Here $V(\xi)$ is a smooth potential that behaves like
$-\frac{1}{4\xi^2}$ as $\xi\to\pm\infty$. It is crucial to notice
that the combined potential behaves like $\frac{2n^2-1/4}{\xi^2}$ as
$\xi\to\pm\infty$. In order to prove our theorems, we express the
resolvent kernel as \[ (\calH_n-(\lambda^2+i0))^{-1}(\xi,\xi')=
\frac{f_+(\xi,\lambda) f_{-}(\xi',\lambda)}{W(\lambda)} \] when
$\xi>\xi'$. Here $f_{\pm}$ are the usual Jost solutions for
$\calH_n$ at energy $\lambda^2$ which are asymptotic to $e^{\pm
i\xi\lambda}$ as $\xi\to\pm\infty$, and $W(\lambda)$ is the
Wronskian of $f_+(\cdot,\lambda)$ with $f_-(\cdot,\lambda)$. It is a
well-known fact of scattering theory, see Deift-Trubowitz~\cite{DT},
that for potentials $V(\xi)$ satisfying $\langle\xi\rangle V(\xi)\in
L^1(\R)$ the Jost solutions exist and are continuous in
$\lambda\in\R$; in fact, they are continuous in $\lambda\ne0$ under
the weaker condition $V\in L^1$. Here, this continuity property --
as well as the existence statement -- fail at $\lambda=0$ because of
the inverse square behavior. Furthermore, it is also common
knowledge that for decay estimates as in Theorem~\ref{thm1}
and~\ref{thm2}, the energy $\lambda=0$ plays a decisive role. For
this reason we need to develop some machinery to determine the
asymptotic behavior of both $f_{\pm}(\cdot,\lambda)$ as well as
$W(\lambda)$ as $\lambda\to0$. We accomplish this by means of two
types of perturbative arguments. The first type is perturbative in
the energy $\lambda$ and around the zero energy solutions; the
latter of course correspond to the harmonic functions on~$\calS$ of
which there are really two for each $n$: if $n=0$ the first one is
constant and the second is logarithmic -- more precisely, it behaves
like $\log\xi$ as $\xi\to\infty$ and like $-\log|\xi|$ as
$\xi\to-\infty$. Starting from these two, we build a fundamental
system of solutions to $\calH_0 f=\lambda^2 f$ at least in the range
$|\xi|\ll \lambda^{-1}$. The second type uses the operator
\[ \widetilde\calH_0:= -\partial_\xi^2 - \frac{1}{4\xi^2} \]
as approximating operator. The Jost solutions of $\widetilde\calH_0$
are given explicitly in terms of (weighted) Hankel functions of
order zero. Using these as approximation, we obtain expressions for
the true Jost solutions which are sufficiently accurate in the range
$|\xi|\gg \lambda^{-\epsilon}$. It is important that this range
overlaps with the range from the previous perturbative argument.
Hence, we are able to glue our fundamental systems together to yield
global Jost solutions. Finally, proving Theorem~\ref{thm1} then
reduces to certain oscillatory integrals for which we rely on
stationary phase type arguments, see~\eqref{4} below. Note that
although these oscillatory integrals are one-dimensional, they still
yield the $t^{-1}$ decay due to the fact that they contain weights
of the form $(\la \xi\ra\la \xi'\ra)^{-\frac12}$.

\section{The proof of Theorem~\ref{thm1}}
\label{sec:proof}

\nind The Laplace-Beltrami operator on $\calS$ is

\beeq\label{1}
\Delta_{\calS}=\frac{1}{r(x)\sqrt{1+r'(x)^2}}\partial_x\left(\frac{r(x)}{\sqrt{1+r'(x)^2}}
\partial_x\right)+\frac{1}{r^2(x)}\partial_{\theta}^2
\eneq
It is convenient to switch to arclength parametrization. Thus,
let

$$\xi (x)=\int_0^x\sqrt{1+r'(y)^2}\, dy.$$

\noindent Then (\ref{1}) can be written as

\beeq\label{2} \Delta_{\calS}=\frac{1}{r(\xi)}\partial_{\xi}(r(\xi
)\partial_{\xi})+\frac{1}{r^2(\xi )}\partial^2_{\theta} \eneq

\noindent where we have abused notation: $r(\xi )$ instead of
$r(x(\xi ))$. We remark that using \eqref{2}, one obtains two
$\theta$ independent harmonic functions on $\calS$:

\begin{align}\label{6}
y_0(\xi )&=1,\nonumber\\
y_1(\xi )&=\int_0^{\xi}r^{-1}(\xi ')\, d\xi '.
\end{align}

\noindent By our asymptotic assumption on $r(x)$,

\beeq\label{7} \xi (x)=\left\{
\begin{array}{ll}
\sqrt{2}\,x+c_{\infty}+O(x^{-1}) & \textrm{as $x\to\infty$}\\
\sqrt{2}\,x-c_{\infty}+O(x^{-1}) & \textrm{as $x\to -\infty$}
\end{array}\right.
\eneq

\noindent where $c_{\infty}$ is some constant.

In particular,

\beeq\label{8}  r(\xi )=\frac{1}{\sqrt{2}\,}\xi
\left(1-\frac{c_{\infty}}{\xi}+O(\xi^{-2})\right)\qquad
\mathrm{as}\; \xi\to\infty.\eneq

Hence,

\beeq\label{9} y_1(\xi )=\sqrt{2}\,\, \sign (\xi )(\log |\xi
|+O(1))\qquad\mathrm{as}\; |\xi |\to\infty. \eneq

\noindent The appearance of $y_1(\xi)$ is remarkable, as it
establishes the existence of a harmonic function on $\mathcal{S}$,
which grows in absolute value like $\log |\xi|$ at both ends. Note
that there can be no harmonic function on~$\R^2$ which grows like
$\log r$ as $r=|\xi|\to\infty$; indeed, this would violate the mean
value property of harmonic functions (on the other hand, the
fundamental solution is $\log r$). In the case of $y_1(\xi)$ it is
crucial that it behaves like $-\log|\xi|$ on one end, and like
$\log|\xi|$ on the other (note the presence of $\sign(\xi)$) --
therefore, $y_1$ is no contradiction to the mean-value property. One
can also think of $\mathcal{S}$ as two planes joined by a neck. Then
on the upper plane the harmonic function grows like $\log|\xi|$,
whereas on the lower plane it
behaves like $-\log|\xi|$.\\

Setting $\omega (\xi ):=\frac{\dot{r}(\xi )}{r(\xi )}$ yields

\beeq\label{2.1}\Delta_{\calS}\,y(\xi ,\theta )=
\partial^2_{\xi}y+\omega\partial_{\xi}y+\frac{1}{r^2}\partial^2_{\theta}y.\eneq

\noindent First we need a few simplifications in order to deal with
the Laplacean. To this end, we remove the first order term in
\eqref{2.1} by setting

\beeq\label{3} y(\xi ,\theta )=r(\xi )^{-1/2}u(\xi ,\theta ). \eneq

\noindent Then

\beeq\label{4}
\partial^2_{\xi}y+\omega\partial_{\xi}y+\frac{1}{r^2}\partial^2_{\theta}y=r^{-1/2}[\partial^2_{\xi}u-V(\xi)
 u+\frac{1}{r^2}\partial^2_{\theta}u]
\eneq

\noindent with

\beeq\label{5} V(\xi )=\frac{1}{4}\omega^2(\xi
)+\frac{1}{2}\dot{\omega}(\xi ). \eneq

Letting $H=-\partial^2_{\xi}+V$, and using \eqref{3} and \eqref{4},
we observe that (\ref{10}) is equivalent with

\beeq\label{11} ||r^{-1/2}e^{itH}r^{-1/2}u||_{L^{\infty}(\R )}\les
t^{-1}||u||_{L^1(\R )}. \eneq

\noindent By the usual reduction to the resolvent of $H$, see Artbazar-Yajima\cite{AY},
Weder\cite{wed}, and Goldberg-Schlag\cite{golsch}, the bound
(\ref{11}) is equivalent to the following oscillatory integral
bound:

\beeq
\begin{aligned}\label{12}
\sup_{\xi >\xi '}&\bigg|\int_0^{\infty}e^{it\lambda^2}(\langle
\xi\rangle \langle \xi'\rangle)^{\frac{-1}{2}}
\lambda\, \text{Im}\left[\frac{f_+(\xi,\lambda)f_-(\xi',\lambda )}{W(\lambda )}\right]\, d\lambda \bigg|\\
&+\sup_{\xi <\xi '}\bigg|\int_0^{\infty}e^{it\lambda^2}(\langle
\xi\rangle \langle \xi'\rangle)^{\frac{-1}{2}}\lambda\,
\text{Im}\left[\frac{f_+(\xi',\lambda)f_-(\xi,\lambda )}{W(\lambda
)}\right]\, d\lambda \bigg|\les t^{-1}
\end{aligned}
\eneq

\noindent where $W(\lambda ):=W(f_+(\cdot ,\lambda),f_-(\cdot
,\lambda))$ is the Wronskian of solutions $f_{\pm}(\cdot ,\lambda )$
of the following ODE

\beeq
\begin{aligned}\label{13}
Hf_{\pm}(\xi ,\lambda )&=-f_{\pm}''(\xi ,\lambda )+V(\xi )f_{\pm}(\xi ,\lambda )\\
&=\lambda^2f_{\pm}(\xi ,\lambda)\\
f_{\pm}(\xi ,\lambda ) &\sim e^{\pm i\lambda\xi}\qquad \mathrm{as}\;
\xi\to\pm\infty
\end{aligned}
\eneq

\noindent provided $\lambda\neq 0$. $f_{\pm}$ are called the {\em
Jost solutions} and it is a standard fact that these solutions
exist. Indeed, since $\omega=\frac{\dot{r}}{r}$, in view of
(\ref{5}) and (\ref{8}), we see that
$$|V(\xi )|\les\langle\xi\rangle^{-2}.$$

In fact, $V$ decays no faster since

\beeq\label{14} V(\xi )=-\frac{1}{4\xi^2}+O(\xi^{-3})\qquad
\mathrm{as}\; |\xi |\to\infty. \eneq It is also important to note
that the term $O(\xi^{-3})$ behave
 like a symbol, i.e. $\vert \frac{d^{k}}{d\xi ^{k}}
 O(\xi^{-3})\vert\lesssim \langle \xi \rangle ^{-3-k}$.
\noindent Therefore, $f_{\pm}(\cdot ,\lambda)$ are solutions of the
{\em Volterra integral equations}

\beeq\label{15} f_+(\xi ,\lambda
)=e^{i\lambda\xi}+\int_{\xi}^{\infty}\frac{\sin(\lambda (\eta
-\xi))}{\lambda}V(\eta )f_+(\eta ,\lambda )\, d\eta \eneq

\noindent and similarly for $f_-$.

However, these integral equations have no meaning at $\lambda =0$.
In fact, the zero energy solutions of $Hu=0$ are given by

\begin{align}\label{16}
u_0(\xi )&=r^{1/2}(\xi ),\nonumber\\
u_1(\xi )&=r^{1/2}(\xi )\int_0^{\xi}r^{-1}(\eta)\, d\eta,
\end{align}

\noindent see (\ref{6}) and (\ref{3}).  Since these functions are
not asymptotically constant as $|\xi |\to +\infty$, it follows that
(\ref{13}) and (\ref{15}) have no meaning at $\lambda=0$.

In passing, we remark that due to the special form of $V$, the
Schr\"odinger operator $H$ can be factorized as

\begin{align}\label{17}
H&=\calL^*\calL,\nonumber\\
\calL&=\frac{d}{d\xi}-\frac{1}{2}\omega.
\end{align}

In particular, $H$ has no negative spectrum.  One can of course
recover $u_0$, $u_1$ in (\ref{16}) by means of (\ref{17}): first,
solve $\calL u_0=0$ and then observe that $\calL^*
(\frac{1}{{u_{0}}})=0$. Therefore, $Hu_0 =0$ and solving $\calL
u_{1} =\frac{1}{u_{0}}$ yields $Hu_1 =0$.

\begin{lem}\label{lemma1}
For any $\lambda\in\R$, define \beeq\label{18} u_j(\xi ,\lambda
):=u_j(\xi
)+\lambda^2\int_0^{\xi}[u_1(\xi)u_0(\eta)-u_1(\eta)u_0(\xi)]u_j(\eta
,\lambda )\, d\eta \eneq \noindent where $j=0,1$.  Then $Hu_j(\cdot
,\lambda )=\lambda^2u_j(\cdot ,\lambda )$ with $u_j(\cdot
,0)=u_j(\cdot )$, for $j=0,1$ and \beeq\label{19} W(u_0(\cdot
,\lambda ),u_1(\cdot ,\lambda ))=1 \eneq for all $\lambda$.
\end{lem}

\begin{proof}
First, one checks that $W(u_0(\cdot ),u_1(\cdot ))=1$.  This yields
$Hu_j(\cdot ,\lambda )=\lambda^2u_j(\cdot ,\lambda )$ since
$Hu_j(\cdot )=0$ for $j=0,1$.  Second, $u_j(0,\lambda )=u_j(0)$ and
$u_j'(0,\lambda )=u_j'(0)$ for $j=0,1$.  Hence $W(u_0(\cdot ,\lambda
),u_1(\cdot ,\lambda ))=u_1'(0)u_0(0)-u_1(0)u_0'(0)=1$.
\end{proof}

As an immediate corollary we have the following statement.

\begin{cor}\label{cor1}
With $f_{\pm}(\cdot ,\lambda )$ as in $(\ref{13})$, one has for any
$\lambda\neq 0$
\begin{align}\label{20}
f_+(\xi ,\lambda )&=a_+(\lambda )u_0(\xi ,\lambda )+b_+(\lambda )u_1(\xi ,\lambda )\nonumber\\
f_-(\xi ,\lambda )&=a_-(\lambda )u_0(\xi ,\lambda )+b_-(\lambda
)u_1(\xi ,\lambda )
\end{align}
\noindent where $a_{\pm}(\lambda )=W(f_{\pm}(\cdot ,\lambda
),u_1(\cdot ,\lambda ))$ and $b_{\pm}(\lambda )=-W(f_{\pm}(\cdot
,\lambda ),u_0(\cdot ,\lambda ))$.  Moreover, $a_-(\lambda
)=a_+(\lambda )$ and $b_-(\lambda )=-b_+(\lambda )$.
\end{cor}

\begin{proof}
The Wronskian relations for $a_{\pm}$, $b_{\pm}$ follow immediately
from (\ref{19}).  Recall that we are assuming $r(x)=r(-x)$ and
therefore also $r(\xi )=r(-\xi )$.  In particular, this implies that
$f_-(-\xi ,\lambda )=f_+(\xi ,\lambda )$ and $u_0(-\xi )=u_0(\xi )$
as well as $u_1(-\xi )=-u_1(\xi )$.  Thus,
$$ a_-(\lambda )= W(f_{-}(\cdot, \lambda),u_1(\cdot ,\lambda))=-W(f_-(-\cdot ,\lambda),u_1(-\cdot ,\lambda))=W(f_+(\cdot ,\lambda ),u_1(\cdot ,\lambda ))=a_+(\lambda )$$
\noindent and
$$ b_-(\lambda )=-W(f_-(\cdot ,\lambda),u_0(\cdot ,\lambda))=
W(f_-(-\cdot ,\lambda),u_0(-\cdot ,\lambda))=W(f_+(\cdot ,\lambda
),u_0(\cdot ,\lambda ))=-b_+(\lambda )$$ as claimed.
\end{proof}
\noindent The point of this corollary is as follows: in order to
prove (\ref{12}) we need to obtain detailed understanding of the
functions $f_{\pm}(\xi ,\lambda )$.  For large $\xi$, we will obtain
asymptotic estimates by perturbing off the potential
$-\frac{1}{4\xi^2}$.  But this analysis fails for small $\xi$, so to
tackle that problem we can use (\ref{18}) to derive useful bounds.
In the end, we need to ``glue'' these two regions together. This is
the meaning of Corollary \ref{cor1}.

\noindent To achieve our goal, namely to show \eqref{12}, we need a
careful analysis of $f_{\pm}(\xi ,\lambda )$. We start by rewriting
(\ref{14}) as

\beeq\label{14'} V(\xi )=-\frac{1}{4\xi^2}+V_1(\xi ),\qquad |\xi |>1
\eneq where $|V_1(\xi )|\les |\xi |^{-3}$ by
Definition~\ref{def:cones}. Moreover, $\vert V_1^{(k)}(\xi
)\vert\les |\xi |^{-3-k}$ for $|\xi |>1$.  Let
$H_0=-\partial^2_{\xi}-\frac{1}{4\xi^2}$.

\begin{lem}\label{lemma4}
For any $\lambda >0$ the problem
\begin{align*}
H_0f_0(\cdot ,\lambda )&=\lambda^2f_0(\cdot ,\lambda ),\\
f_0(\xi ,\lambda )&\sim e^{i\xi\lambda}
\end{align*}
as $\xi\to\infty$ has a unique solution on $\xi >0$.  It is given by
\beeq\label{23} f_0(\xi ,\lambda )=\sqrt{\frac{\pi}{2}}\, e^{i\pi
/4}\sqrt{\xi\lambda}\,H^{(+)}_0(\xi\lambda). \eneq Here
$H^{(+)}_0(z)=J_0(z)+iY_0(z)$ is the Hankel function.
\end{lem}

\begin{proof}
It is well-known, see Abramowitz-Stegun\cite{AS}, that the ODE
$$w''(z)+\left(\lambda^2+\frac{1}{4z^2}\right)W(z)=0$$
has a fundamental system of solutions $\sqrt{z}\, J_0(\lambda z)$, $\sqrt{z}\, Y_0(\lambda z)$ or equivalently, $\sqrt{z}\, H^{(+)}_0(\lambda z)$, $\sqrt{z}\, H^{(-)}_0(\lambda z)$.\\

Recall the asymptotics

\begin{align*}
H^{(+)}_0(x)&\sim\sqrt{\frac{2}{\pi x}}\, e^{i(x-\frac{\pi}{4})}
\qquad\mathrm{as}\; x\to +\infty\\
H^{(-)}_0(x)&\sim\sqrt{\frac{2}{\pi x}}\,
e^{-i(x-\frac{\pi}{4})}\qquad\mathrm{as}\; x\to +\infty.
\end{align*}

Thus, (\ref{23}) is the unique solution so that $$f_0(\xi ,\lambda
)\sim e^{i\xi\lambda},$$ as claimed.
\end{proof}

As was observed in \eqref{15}, The Volterra type integral equations
play a crucial role in our analysis throughout the paper so for the
convenience of the reader we will very briefly sketch how one solves
such equations. Let us consider the following Volterra equations
\begin{equation*}
(\ast )\,\,\,f(x)= g(x)+\int_{x}^{\infty} K(x,s) f(s) ds,
\end{equation*}
or
\begin{equation*}
(\ast\ast )\,\,\,f(x)= g(x)+\int_{a}^{x} K(x,s) f(s) ds,
\end{equation*}
with some $g(x)\in L^{\infty}$ and $a\in \mathbb{R}$. Evidently, one
solves them by an iteration procedure which requires finding a
suitable convergent majorant for the resulting series expansion. In
Lemma~\ref{4'} below we show that, depending on the choice of norm
of $K$, this majorant is either a geometric series or an exponential
series. We remark that -- as usual -- one only needs the latter
alternative in this paper since it does not require a smallness
condition.

\begin{lem}\label{4'}
Let $a\in \mathbb{R}$ and $g(x)\in L^{\infty}(a,\infty)$.

\begin{itemize}
\item
  If\, $M:=\sup_{x>a}\int_{x}^{\infty}\vert K(x,s)\vert ds<1$, then there exists a unique $L^{\infty}$ solution to the
Volterra equation $(\ast)$, valid for $x>a$, which is given by
\begin{equation*}
f(x)= g(x)+\sum_{n=1}^{\infty}\int_{x}^{\infty} K^{n}(x,s)g(s)\,ds,
\end{equation*}
with $K^1(x,s)= K(x,s)$ and $K^{n}(x,s)
=\int_{x}^{\infty}K^{n-1}(x,t)K(t,s) dt$ for $n\geq 2.$ One also has
the bound
\[
\|f\|_{L^\infty(a,\infty)} \le (1-M)^{-1}
\|g\|_{L^\infty(a,\infty)},
\]
and a similar statement holds for $(\ast\ast)$.

\item Let $\mu:= \int_{a}^\infty \sup_{a<x<s} |K(x,s)|\, ds<\infty$. Then
there exists a unique solution to $(\ast)$ given by
\begin{equation}
\label{eq:volt_it} f(x) = g(x) + \sum_{n=1}^\infty \int_a^\infty
\ldots \int_a^\infty \prod_{i=1}^n \chi_{[x_{i-1}<x_i]}
K(x_{i-1},x_i) \; g(x_{n}) \, dx_{n}\ldots dx_1.
\end{equation}
Furthermore one has the bound
\[ \|f\|_{L^\infty(a,\infty)} \le e^\mu \|g\|_{L^\infty(a,\infty)}, \]
and an analogue statement holds for $(\ast\ast)$.
\end{itemize}
\end{lem}
\begin{proof}
We only prove the lemma for $(\ast)$ since the proof for
$(\ast\ast)$ is almost identical. The solutions to both equations
are found through a standard Picard-Volterra iteration procedure,
which we will from now on refer to as the {\em Volterra iteration}.
For $(\ast )$, this iteration yields a solution via the
Picard-Banach fixed point theorem. We also have \beeq\label{25'}
\vert \int_{x}^{\infty} K^{n}(x,s)g(s)ds\vert\leq M^{n} \Vert
g\Vert_{L^{\infty}}, \eneq for $x>a$. So if $M<1$ then the series
representing the solution $f(x)$, will converge absolutely and
uniformly for $x>a$ and the formula for the sum of a geometric
series provides us with the bound for the $L^{\infty}$ norm of $f$.

For the second part, we show that the infinite Volterra iteration
\eqref{eq:volt_it} for $(\ast)$ converges. To this end, define
\[ K_0(s) := \sup_{a<x<s} |K(x,s)| \]
Then
\begin{align*}
& \Big| \int_a^\infty \ldots \int_a^\infty
\prod_{i=1}^n \chi_{[x_{i-1}<x_i]} K(x_{i-1},x_i) \; g(x_{n}) \, dx_{n}\ldots dx_1 \Big| \\
& \le \int_a^\infty \ldots \int_a^\infty
\prod_{i=1}^n \chi_{[x_{i-1}<x_i]} K_0(x_i) \; |g(x_{n})| \, dx_{n}\ldots dx_1  \\
& = \|g\|_{L^\infty(a,\infty)}\frac{1}{n!} \int_a^\infty \ldots \int_a^\infty
\prod_{i=1}^n  K_0(x_i) \, dx_{n}\ldots dx_1  \\
&= \frac{1}{n!} \|g\|_{L^\infty(a,\infty)} \Big(\int_a^\infty K_0(s)\, ds\Big)^n
\end{align*}
Hence, the series in \eqref{eq:volt_it} converges absolutely and uniformly in $x>a$ with
the uniform upper bound
\[ \|g\|_{L^\infty(a,\infty)}\sum_{n=0}^\infty \frac{1}{n!}\mu^n = e^\mu\|g\|_{L^\infty(a,\infty)}\]
as claimed.
\end{proof}

Having these tools at our disposal, we proceed with our
investigation of the Jost solutions. To this end, instead of the
Volterra equation (\ref{15}) we will work with the following
representation of the solutions of (\ref{13}):

\begin{lem}\label{lemma5}
For any $\xi >0$, $\lambda >0$, \beeq\label{24} f_+(\xi ,\lambda
)=f_0(\xi ,\lambda )+\int_{\xi}^{\infty}G_0(\xi ,\eta ;\lambda
)V_1(\eta )f_+(\lambda,\eta)\, d\eta \eneq with $V_1$ as in
$(\ref{14'})$, $f_0$ as in $(\ref{23})$ and \beeq\label{25} G_0(\xi
,\eta ;\lambda)=[\overline{f_0(\xi ,\lambda )}f_0(\lambda,\eta
)-f_0(\xi ,\lambda)\overline{f_0(\lambda,\eta )}](-2i\lambda )^{-1}.
\eneq
\end{lem}

\begin{proof}
Simply observe that $G_0$ is the Green's function of our problem
relative to $H_0$.  Indeed,
\begin{align*}
G_0(\xi ,\xi;\lambda )&=0,\\
\partial_{\xi}G_0(\xi,\eta;\lambda)|_{\eta =\xi}&=1,\\
H_0G_0(\cdot,\eta;\lambda )&=\lambda^2G_0(\cdot,\eta;\lambda).
\end{align*}

Here we have used that $W(f_0(\cdot,\lambda ),\overline{f_0(\cdot,\lambda)})=-2i\lambda$ which can be seen by computing the Wronskian at $\xi =\infty$.\\

In conclusion,
$$H_0f_+(\xi,\lambda )=\lambda^2\left[f_0(\xi,\lambda)+\int_{\xi}^{\infty}G_0(\xi,\eta;\lambda)V_1(\eta)f_+(\lambda,\eta)\, d\eta \right]-V_1(\xi )f_+(\xi,\lambda )$$
or equivalently,
$$Hf_+(\cdot,\lambda )=\lambda^2f_+(\cdot ,\lambda ).$$

Finally, observe that for $\xi >\lambda^{-1}$ fixed,
$$\sup_{\eta >\xi}|G_0(\xi,\eta ;\lambda )|\les\lambda^{-1}.$$
By the Volterra iteration discussed above, this implies that
$|f_+(\xi,\lambda )-f_0(\xi,\lambda )|\les \lambda^{-1}\xi^{-2}.$ In
particular,
$$f_+(\xi,\lambda )\sim e^{i\lambda\xi} \qquad\mathrm{as}\; \xi\to\infty$$
and we are done.
\end{proof}

For small arguments the Hankel function $H_0(z)$ displays the
following asymptotic behavior, see \cite{AS}:

\beeq\label{26} H^{(+)}_0(z)=1+O_{\R}(z^2)+\frac{2}{\pi}i\log
z+i\varkappa +iO_{\R}(z^2\log z) \eneq

\noindent as $z\to 0$ where $\varkappa$ is some real constant.\\

\noindent Estimating the oscillatory integrals will require
understanding
$\partial^{k}_{\lambda}\partial^{l}_{\xi}f_{\pm}(\xi,\lambda)$, for
$0\leq k+l\leq 2$, $W(\lambda)$, $W'(\lambda )$ and thus
$a_{\pm}(\lambda)$, $b_{\pm}(\lambda)$, $a_{\pm}'(\lambda )$ and
$b_{\pm}'(\lambda )$. To obtain asymptotic expansions for all these
functions, we need to know the asymptotic behavior of $u_{j}(\xi)$,
and thereafter that of
$\partial^{k}_{\lambda}\partial^{l}_{\xi}u_{j}(\xi,\lambda)$, for
$j=1,\, 2$ and $0\leq k+l\leq 2$.

We start by analyzing the $u_{j}(\xi)$'s.
\begin{lem}\label{lemma7} \beeq\label{29} r(\xi
)=\frac{1}{\sqrt{2}\,}\xi
\left(1-\frac{c_{\infty}}{\xi}+O(\xi^{-2})\right)\qquad
\mathrm{as}\; \xi\to\infty. \eneq With $c_{\infty}$ as in
$\eqref{7}$, In particular,
\begin{align}\label{30}
u_0(\xi )&=2^{-1/4}\xi^{1/2}\left(1-\frac{c_{\infty}}{2\xi}+O(\xi^{-2})\right)\qquad \mathrm{as}\; \xi\to\infty\nonumber\\
u_1(\xi
)&=2^{1/4}\xi^{1/2}\left(1-\frac{c_{\infty}}{2\xi}+O(\xi^{-2})\right)\left(\log\xi+c_2+O(\xi^{-1})\right).
\end{align}
Here $c_2$ is some constant. Moreover, the $O$-terms behave like
symbols.
\end{lem}

\begin{proof}
In view of \eqref{0} and \eqref{7},
\[
\begin{array}{l}
\xi (x)=\sqrt{2}\,x+c_{\infty}+O(x^{-1})\\
r(x)=x+O(x^{-1})\end{array} \textrm{both as } x\to\infty
\]
where the $O$-terms behave like symbols.  This implies (\ref{29}), as well as the expansion of $u_0(\xi )=\sqrt{r(\xi )}$ in (\ref{30}).\\
Next compute
\begin{align*}
\int_0^{\xi}r^{-1}(\eta )\, d\eta &= \int_0^{\xi}\sqrt{2}\,\langle\eta\rangle^{-1}\left(1+\frac{c_{\infty}}{\langle\eta\rangle}
+O(\langle\eta\rangle^{-2})\right)\, d\eta\\
&=\sqrt{2}\,(\log\xi +c_2)+O(\xi^{-1})\qquad\mathrm{as}\;
\xi\to\infty.
\end{align*}
Thus,
\begin{align*}
u_1(\xi )&=\sqrt{r(\xi )}\int_0^{\xi}r^{-1}(\eta )\, d\eta\\
&=2^{1/4}\xi^{1/2}\left(1-\frac{c_{\infty}}{2\xi}+O(\xi^{-2})\right)\left(\log\xi+c_2+O(\xi^{-1})\right)\qquad\mathrm{as}\;
\xi\to\infty.
\end{align*}
\end{proof}

To study the behavior of $u_{j}(\xi ,\lambda)$'s we recall the
Volterra equation \eqref{18}
\begin{equation*}
u_j(\xi ,\lambda ):=u_0(\xi
)+\lambda^2\int_0^{\xi}[u_1(\xi)u_0(\eta)-u_1(\eta)u_0(\xi)]u_j(\eta
,\lambda )\, d\eta.
\end{equation*}

Hence, setting $h_j(\xi,\lambda):=\frac{u_j(\xi ,\lambda )}{u_j(\xi
)}$, for $\xi >0$ we obtain the integral equations

\beeq \label{31'} h_{0}(\xi ,\lambda
)=1+\frac{\lambda^2}{u_{0}(\xi)}\int_0^{\xi}[u_1(\xi)u^{2}_0(\eta)-u_0(\xi)u_1(\eta)u_0
(\eta)]h_0(\eta ,\lambda )\, d\eta, \eneq

\beeq \label{31''} h_{1}(\xi ,\lambda
)=1+\frac{\lambda^2}{u_{1}(\xi)}\int_0^{\xi}[u_1(\xi
)u_0(\eta)u_1(\eta )-u_0(\xi )u_1^2(\eta )] h_{1}(\eta ,\lambda)\,
d\eta. \eneq

Therefore, in order to solve the integral equation for
$u_{j}(\xi,\lambda)$ for large $\xi$, it is enough to carry out the
Volterra iteration for \eqref{31'} and \eqref{31''} which are
simpler since the first iterates in both cases are identically equal
to $1$, and then multiply the $h_j$'s so obtained by the
$u_{j}(\xi)$'s of \eqref{30}. For the Volterra iterations we would
need to understand the behavior of the kernels in the integral
equations \eqref{31'} and \eqref{31''}. For this purpose, Lemma
\ref{lemma7} yields

\begin{cor}\label{cor3}
As $\xi\to\infty$,
\begin{align}\label{31}
u_1(\xi )\int_0^{\xi}u_0^2(\eta )\, d\eta -u_0(\xi
)\int_0^{\xi}u_1u_0(\eta )\, d\eta &=
\frac{1}{4}2^{-1/4}\xi^{5/2}+O(\xi^{3/2}\log\xi )
\end{align}
\begin{align}\label{32}
u_1(\xi )\int_0^{\xi}u_0u_1(\eta )\, d\eta -u_0(\xi
)\int_0^{\xi}u_1^2(\eta )\, d\eta &=
\frac{1}{4}2^{1/4}\xi^{5/2}\log\xi
+c_3\xi^{5/2}+O(\xi^{\frac{3}{2}}\log\xi )
\end{align}
\end{cor}

\begin{proof}
In view of the first equality of (\ref{30}) and a justifiable
modification of that expression at $0$, we have
\begin{align*}
\int_0^{\xi}u_0^2(\eta )\, d\eta & = 2^{-1/2}\int_0^{\xi}\eta \left(1-\frac{c_{\infty}}{\langle\eta\rangle}+O(\langle\eta\rangle^{-2})\right)\,d\eta\\
&=2^{-1/2}\left(\frac{1}{2}\xi^2-c_{\infty}\cdot\xi + O(\log\xi )\right)\\
\int_0^{\xi}u_0(\eta )u_1(\eta )\, d\eta & = \int_0^{\xi}\eta \left(1-\frac{c_{\infty}}{\langle\eta\rangle}+O(\langle\eta\rangle^{-2})\right)\left(\log\eta +c_2+O(\langle\eta\rangle^{-1})\right)\,d\eta\\
&=\frac{1}{2}\xi^2\log\xi
+\frac{1}{2}\left(c_2-\frac{1}{2}\right)\xi^2+O(\xi\log\xi ).
\end{align*}
Thus,
\begin{align*}
(\ref{31})=&2^{-1/4}\xi^{1/2}(\log\xi +c_2+O(\xi^{-1}\log\xi ))\left(\frac{1}{2}\xi^2+O(\xi )\right)\\
&-2^{-1/4}\xi^{1/2}(1+O(\xi^{-1}))\left(\frac{1}{2}\xi^2\log\xi+\frac{1}{2}\left(c_2-\frac{1}{2}\right)\xi^2+O(\xi\log\xi )\right)\\
=&2^{-1/4}\xi^{1/2}\left[\frac{1}{4}\xi^2+O(\xi\log\xi )\right]
\end{align*}
Next, compute
\begin{align*}
\int_0^{\xi}u_1^2(\eta )\, d\eta &=\sqrt{2}\,\int_0^{\xi}\eta (\log^2\eta +2c_2\log\eta +O(\langle\eta\rangle^{-1}\log\eta ))(1+O(\langle\eta\rangle^{-1}))\, d\eta\\
&=\sqrt{2}\,\left(\frac{1}{2}\xi^2\log^2\xi +(2c_2-1)\int_0^{\xi}\eta\log\eta\, d\eta +O(\xi\log^2\xi )\right)\\
&=\sqrt{2}\,\left(\frac{1}{2}\xi^2\log^2\xi
+\frac{2c_2-1}{2}\xi^2\log\xi -\frac{2c_2-1}{4}\xi^2+O(\xi\log^2\xi
)\right)
\end{align*}
Thus, $(\ref{32})=$
\begin{align*}
2^{1/4}&\xi^{1/2}(\log\xi +c_2+O(\xi^{-1}))(1+O(\xi^{-1}))
\left(\frac{1}{2}\xi^2\log\xi
+\frac{1}{2}\left(c_2-\frac{1}{2}\right)\xi^2
+O(\xi\log\xi )\right)\\
&-2^{1/4}\xi^{1/2}(1+O(\xi^{-1})) \left(\frac{1}{2}\xi^2\log^2\xi
+\frac{2c_2-1}{2}\xi^2\log\xi -\frac{2c_2-1}{4}\xi^2+O(\xi\log^2\xi )\right)\\
=2^{1/4}&\xi^{1/2}\bigg\{\frac{1}{2}\xi^2\log^2\xi +\frac{2c_2-\frac{1}{2}}{2}\xi^2\log\xi +O(\xi\log^2\xi )+\frac{c_2}{2}\left( c_2-\frac{1}{2}\right)\xi^2\\
&-\frac{1}{2}\xi^2\log^2\xi -\frac{2c_2-1}{2}\xi^2\log\xi +\frac{2c_2-1}{4}\xi^2\bigg\}\\
=2^{1/4}&\sqrt{\xi}\left(\frac{1}{4}\xi^2\log\xi+2^{-1/4}c_3\xi^2+O(\xi\log\xi
)\right)
\end{align*}
as claimed.
\end{proof}

Thus a Volterra iteration and the preceding yields the following
result for the $u_{j}(\xi,\lambda)$'s. The importance of
Corollary~\ref{cor4} lies with the fact that we do not lose
$\log\xi$ factors in the $O(\cdot)$-terms
 as such factors would destroy the dipsersive estimate. It is easy to see that carrying
out the Volterra iteration crudely, by putting absolute values inside the integrals, leads
to such $\log\xi$ losses. Therefore, we actually need to compute the Volterra iterates in~\eqref{eq:volt_it} explicitly
(or, more precisely, its analogue for ($**$))
.

\begin{cor}\label{cor4}
In the range $1<<\xi <<\lambda^{-1}$, $j=0,1$,
\begin{align}\label{33}
u_j(\xi ,\lambda )&=u_j(\xi )(1+O((\xi\lambda )^2))\\
\partial_{\xi}u_j(\xi ,\lambda )&=u_j'(\xi )(1+O((\xi\lambda )^2))\nonumber
\end{align}
\begin{align}\label{34}
\partial_{\lambda}u_0(\xi ,\lambda )&=\frac{1}{2}2^{-1/4}\lambda (\xi^{5/2}+O(\xi^{3/2}\log\xi ))(1+O((\xi\lambda )^2))\\
\partial_{\lambda}u_1(\xi ,\lambda )&=\frac{1}{2}2^{1/4}\lambda (\xi^{5/2}\log\xi +c_3\xi^{5/2}+O(\xi^{3/2}\log\xi ))(1+O((\xi\lambda )^2))\nonumber
\end{align}
\begin{align}\label{35}
\partial^2_{\lambda\xi}u_0(\xi ,\lambda )&=\frac{5}{4}2^{-1/4}\lambda (\xi^{3/2}+O(\xi^{1/2}\log\xi ))(1+O((\xi\lambda )^2))\\
\partial^2_{\lambda\xi}u_1(\xi ,\lambda )&=\frac{5}{4}2^{1/4}\lambda (\xi^{3/2}\log\xi +\frac{2}{5}\xi^{3/2}+c_3\xi^{3/2}+O(\xi^{1/2}\log\xi ))(1+O((\xi\lambda )^2))\nonumber
\end{align}
\end{cor}

\begin{proof}
We sketch the proof of this somewhat computational lemma, for the
function $u_{1}(\xi,\lambda)$ since the argument for
$u_{0}(\xi,\lambda)$ is completely anaolgous and in fact easier. The proof of the first
equality in \eqref{33} is based on the Volterra integral equation
\eqref{31''} \beeq h_{1}(\xi ,\lambda
)=1+\lambda^2\int_0^{\xi}[\frac{u_1(\xi )u_0(\eta)u_1(\eta )-u_0(\xi
)u_1^2(\eta )}{u_{1}(\xi)}] h_{1}(\eta ,\lambda)\, d\eta \eneq and
its derivatives in both $\xi$ and $\lambda$ and the Volterra
iteration, for which we also need to use Corollary \ref{cor3}. The
iteration will produce a solution which is given by
\begin{align*}
 h_{1}(\xi ,\lambda) = & 1 + \sum_{n=1}^\infty \lambda^{2n} \int_0^\xi
\int_{0}^{\xi_{1}}\ldots \int_{0}^{\xi_{n-1}} \frac{u_1(\xi
)u_0(\xi_{1})u_1(\xi_{1} )-u_0(\xi )u_1^2(\xi_{1} )}{u_{1}(\xi)}\cdots\\
&\frac{u_1(\xi_{n-1} )u_0(\xi_{n})u_1(\xi_{n} )-u_0(\xi_{n-1}
)u_1^2(\xi_{n})}{u_{1}(\xi_{n-1} )} \, d\,\xi_{n}\ldots d\,\xi_1 = \\
& 1+ \lambda^{2}\int_0^\xi  \frac{u_1(\xi )u_0(\xi_{1})u_1(\xi_{1}
)-u_0(\xi )u_1^2(\xi_{1} )}{u_{1}(\xi)} d\,\xi_{1}+\\
&\lambda^{4}\int_{0}^{\xi}\int_{0}^{\xi_{1}} \frac{u_1(\xi
)u_0(\xi_{1})u_1(\xi_{1} )-u_0(\xi )u_1^2(\xi_{1}
)}{u_{1}(\xi)}\times\\&\frac{u_1(\xi_{1} ) u_0(\xi_{2})u_1(\xi_{2}
)-u_0(\xi_{1} )u_1^2(\xi_{2})}{u_{1}(\xi_{1} )} \, d\,\xi_{2}\,
d\,\xi_1 +\cdots
\end{align*}
Therefore, \eqref{32} and the equalities
\begin{align*}
u_0(\xi )&=2^{-1/4}\xi^{1/2}\left(1-\frac{c_{\infty}}{2\xi}+O(\xi^{-2})\right)\nonumber\\
u_1(\xi
)&=2^{1/4}\xi^{1/2}\left(1-\frac{c_{\infty}}{2\xi}+O(\xi^{-2})\right)\left(\log\xi+c_2+O(\xi^{-1})\right)
\end{align*}
yield
\begin{align*}
h_{1}(\xi ,\lambda)= & 1+
\frac{\lambda^{2}}{u_{1}(\xi)}(\frac{1}{4}2^{1/4}\xi^{5/2}\log\xi
+c_3\xi^{5/2}+O(\xi^{\frac{3}{2}}\log\xi ))+\\
&\lambda^{4}\{\int_{0}^{\xi}
u_0(\xi_{1})[\frac{1}{4}2^{1/4}\xi_{1}^{5/2}\log\xi_{1}
+c_3\xi_{1}^{5/2}+O(\xi_{1}^{\frac{3}{2}}\log\xi_{1})] d\xi_{1}-\\
&\frac{u_{0}(\xi)}{u_{1}(\xi)}\int_{0}^{\xi}u_1(\xi_{1}
)[\frac{1}{4}2^{1/4}\xi_{1}^{5/2}\log\xi_{1}
+c_3\xi_{1}^{5/2}+O(\xi_{1}^{\frac{3}{2}}\log\xi_{1})]\,d\xi_{1}\}
+\cdots = 1+ O(\lambda^{2}\xi^{2}),
\end{align*}
since we are assuming that $1<<\xi <<\lambda^{-1}$.
The point to notice here is that terms involving $\xi^4\log\xi$ (the leading orders)
after the integration cancel.
Furthermore,
we obtain the usual $n!$ gain from the Volterra
iteration, see Lemma~\ref{4'}, from repeated integration of powers rather
than from symmetry considerations.
Hence
$u_{1}(\xi,\lambda)= u_{1}(\xi)(1+O(\lambda^{2}\xi^{2}))$ in that
range. To deal with the derivatives, it is more convenient to
differentiate directly the integral equation \eqref{18} for
$u_{1}(\xi,\lambda)$ with respect to $\xi$ and/or $\lambda$, which yields respectively
\begin{align}\label{iter1}
\partial_{\xi}u_{1}(\xi,\lambda)= \partial_{\xi} u_{1}(\xi) + \lambda ^2 \int_{0}^{\xi}
[\partial_{\xi}u_1 (\xi) u_{0}(\eta)-u_{1}(\eta)\partial_{\xi}
u_{0}(\xi)]u_{1}( \eta, \lambda) d\eta
\end{align}

\beeq\label{iter2}
\begin{split}
\partial_{\lambda}u_{1}(\xi,\lambda)= & 2\lambda \int_{0}^{\xi}[u_1
(\xi) u_{0}(\eta)-u_{1}(\eta) u_{0}(\xi)]u_{1}( \eta, \lambda) d\eta
+\\ & \lambda^{2}\int_{0}^{\xi}[u_1 (\xi)
u_{0}(\eta)-u_{1}(\eta)u_{0}(\xi)]\partial_{\lambda}u_{1}(\eta,
\lambda) d\eta,
\end{split}
\end{equation}
and \beeq\label{iter3}
\begin{split}
\partial^{2}_{\lambda\, \xi}u_{1}(\xi,\lambda)= & 2\lambda
\int_{0}^{\xi}[\partial_{\xi}u_1 (\xi)
u_{0}(\eta)-u_{1}(\eta)\partial_{\xi} u_{0}(\xi)]u_{1}( \eta,
\lambda) d\eta +\\ & \lambda^{2}\int_{0}^{\xi}[\partial_{\xi}u_1
(\xi)
u_{0}(\eta)-u_{1}(\eta)\partial_{\xi}u_{0}(\xi)]\partial_{\lambda}u_{1}(\eta,
\lambda) d\eta.
\end{split}
\end{equation}
In dealing with \eqref{iter1}, we simply plug in the information
from the first equality of \eqref{33} and calculate the resulting
integral. For \eqref{iter2}, we observe that by \eqref{32} the term
\[ 2\lambda \int_{0}^{\xi}[u_1 (\xi) u_{0}(\eta)-u_{1}(\eta)
u_{0}(\xi)]u_{1}( \eta, \lambda) d\eta\] is equal to $\lambda
(\frac12 2^{1/4}\xi^{5/2}\log\xi +2c_3\xi^{5/2}+O(\xi^{3/2}\log\xi
)) $. Therefore to solve \eqref{iter2}, one needs to run the
Volterra iteration with this expression as the first iterate. The
treatment of \eqref{iter3} is similar to that of \eqref{iter2} and
we skip the details.
\end{proof}

We now turn to $f_{\pm}(\xi,\lambda )$ as well as $a_{\pm}$,
$b_{\pm}(\lambda )$.

\begin{lem}\label{lemma8}
If $\lambda >0$ is small, and $1<<\xi <<\lambda^{-1}$, then
\beeq\label{35'} f_+(\xi,\lambda )=f_0(\xi,\lambda
)+O(\xi^{-1/2}\lambda^{\frac{1}{2}-\eps})\nonumber \eneq with $\eps
>0$ arbitrary.
\end{lem}

\begin{proof}
We observed above that, c.f. (\ref{24}) and (\ref{25}), \beeq
|G_0(\xi,\eta;\lambda )|\les\sqrt{\xi\eta}\,|\log\lambda
|^2\,\chi_{[\xi <\eta
<\lambda^{-1}]}+\sqrt{\frac{\xi}{\lambda}}\,|\log\lambda
|\,\chi_{[\eta >\lambda^{-1}]}\nonumber \eneq Thus integrating and
taking $1<<\xi <<\lambda^{-1}$ into account, we obtain
\begin{align*}
\bigg|\int_{\xi}^{\infty}G_0(\xi,\eta;\lambda )V_1(\eta )f_+(\eta,\lambda )\, d\eta \bigg|\les&\int_{\xi}^{\lambda^{-1}}\sqrt{\xi\eta}\,|\log\lambda |^2\,\eta^{-3}\sqrt{\eta\lambda}\,|\log\lambda |^3\, d\eta\\
&+\int_{\lambda^{-1}}^{\infty}\sqrt{\frac{\xi}{\lambda}}\,|\log\lambda |\,\eta^{-3}\, d\eta\\
\les&\xi^{-1/2}\lambda^{\frac{1}{2}-\eps},
\end{align*}
as claimed.
\end{proof}

\begin{lem}\label{lemma9}
For small $\lambda >0$,
\begin{align}\label{36}
a_+(\lambda )&=2^{1/4}c_0\sqrt{\lambda}(1+ic_1\log\lambda +ic_3)+O(\lambda^{1-\eps})\\
b_+(\lambda
)&=i2^{-1/4}c_0c_1\sqrt{\lambda}+O(\lambda^{1-\eps}),\nonumber
\end{align}
where $c_0=\sqrt{\frac{\pi}{2}}e^{i\frac{\pi}{4}}$,
$c_1=\frac{2}{\pi}$, and $c_3$ is some real constant.
\end{lem}

\begin{proof}
By Corollary \ref{cor1} we have $a_+(\lambda )=f_+(\xi,\lambda
)u'_1(\xi,\lambda )-f'_+(\xi,\lambda )u_1(\xi,\lambda )$. Hence
Lemma \ref{lemma8} and Corollary \ref{cor4} with
$\xi=\lambda^{-1/2}$ yield,
\begin{align*}
c_0^{-1}2^{1/4}a_+=&\sqrt{\lambda\xi}\,H_0(\xi\lambda )\frac{1}{2}\xi^{-1/2}(\log\xi +c_2+2)\\
&-\left(\frac{1}{2}\xi^{-1/2}\sqrt{\lambda}\,H_0(\xi\lambda )+\sqrt{\xi\lambda}\,H_0'(\xi\lambda )
\lambda \right)\xi^{1/2}(\log\xi +c_2)+O(\lambda^{1-\eps})\\
=&\sqrt{\lambda}\,H_0(\xi\lambda )-\sqrt{\xi\lambda}\: \frac{ic_1}{\xi}\sqrt{\xi}(\log\xi +c_2)+O(\lambda^{1-\eps})\\
=&\sqrt{\lambda}(1+ic_1\log (\xi\lambda )+i\varkappa -ic_1\log\xi -ic_1c_2)+O(\lambda^{1-\eps})\\
=&\sqrt{\lambda}(1+ic_1\log\lambda +ic_3)+O(\lambda^{1-\eps}),
\end{align*}
as claimed.  Note that $c_3=\varkappa -c_1c_2$.\\

Similarly,
\begin{align*}
-c_0^{-1}2^{\frac{1}{4}}b_+&=\sqrt{\lambda\xi}\, H_0(\xi\lambda )
\frac{1}{2}\xi^{-1/2}-\xi^{1/2}\left(\frac{1}{2}\xi^{-1/2}
\sqrt{\lambda}\, H_0(\lambda\xi )+\sqrt{\xi\lambda}\, H_0'(\xi\lambda )\lambda \right)\\
&+O(\lambda^{1-\eps})\\
&=-\xi\sqrt{\lambda}\:\frac{ic_1}{\xi\lambda}\,\lambda+O(\lambda^{1-\eps})\\
&=-ic_1\sqrt{\lambda}+O(\lambda^{1-\eps}),
\end{align*}
and the lemma follows.
\end{proof}

Using the expressions for $a_{+}$ and $b_{+}$ above, we obtain the
following

\begin{cor}\label{cor5}
Let $\lambda >0$ be small.  Then \beeq\label{37} f_+(\xi, \lambda
)=c_0\sqrt{\lambda\langle\xi\rangle}\left(1+ic_1\log
(\lambda\langle\xi\rangle
)+ic_4+O(\lambda^{\frac{1}{2}-\eps})+O(\langle\xi\rangle^{-1}\log\langle\xi\rangle
)\right) \eneq for $0<\xi <\lambda^{-1}$, whereas for
$-\lambda^{-1}<\xi <0,$ \beeq\label{38} f_+(\xi, \lambda
)=c_0\sqrt{\lambda\langle\xi\rangle}\left(1+ic_1\log
(\lambda\langle\xi\rangle^{-1}
)+ic_5+O(\lambda^{\frac{1}{2}-\eps})+O(\langle\xi\rangle^{-1}\log\langle\xi\rangle
)\right) \eneq
\end{cor}

\begin{proof}
This follows by inserting our asymptotic expansions for $a_+(\lambda )$, $b_+(\lambda )$,\\
 $u_0(\xi,\lambda )$, and $u_1(\xi,\lambda )$ into (\ref{20}).
\end{proof}

 We also need some information about certain partial derivatives of
 $f_{+}(\xi,\lambda)$. This is provided by

\begin{lem}\label{lemma10}
For $\lambda >0$ small and $1<<\xi <<\lambda^{-1}$ we have
\begin{align*}
\partial_{\xi}f_+(\xi,\lambda )&=\partial_{\xi} f_0(\xi,\lambda )+O(\xi^{-3/2}\lambda^{\frac{1}{2}-\eps})\\
\partial_{\lambda}f_+(\xi,\lambda )&=\partial_{\lambda} f_0(\xi,\lambda )+O(\xi^{-1/2}\lambda^{-\frac{1}{2}-\eps})\\
\partial^2_{\xi\lambda}f_+(\xi,\lambda )&=\partial^2_{\xi\lambda} f_0(\xi,\lambda )+O(\xi^{-3/2}\lambda^{-\frac{1}{2}-\eps})
\end{align*}
with $\eps >0$ arbitrary.
\end{lem}

\begin{proof}
This follows by taking derivatives in Lemma \ref{lemma8}.
\end{proof}

To be able to carry out the analysis, one also needs to understand
the derivative of the Wronskian. To that end we have

\begin{cor}\label{cor6}
For small $\lambda >0$,
\begin{align}\label{39}
a_+'(\lambda )&=\frac{1}{2}2^{1/4}c_0\lambda^{-1/2}(1+ic_3+2ic_1+ic_1\log\lambda )+O(\lambda^{-\eps})\\
b'_+(\lambda
)&=\frac{i}{2}2^{-1/4}c_0c_1\lambda^{-1/2}+O(\lambda^{-\eps})\nonumber
\end{align}
where $\eps >0$ is arbitrary.
\end{cor}

\begin{proof}
In view of the preceding,
\begin{align}\label{40}
a'_+(\lambda )=&W(\partial_{\lambda}f_+,u_1)+W(f_+,\partial_{\lambda}u_1)\nonumber\\
=&W(\partial_{\lambda}f_0,u_1)+W(f_0,\partial_{\lambda}u_1)+O(\lambda^{-\eps})\nonumber\\
=&\partial_{\lambda}[c_0\sqrt{\lambda\xi}H_0(\lambda\xi )]\left(\frac{1}{2}\xi^{-1/2}(\log\xi +c_2)+\xi^{-1/2}\right)2^{1/4}\\
&-\partial^2_{\lambda\xi}[c_0\sqrt{\lambda\xi}H_0(\lambda\xi )]\xi^{1/2}(\log\xi +c_2)\cdot 2^{1/4}\nonumber\\
&+c_0\sqrt{\lambda\xi}H_0(\lambda\xi )\cdot\frac{5}{4}\cdot 2^{1/4}\lambda \left(\xi^{3/2}\log\xi +\left(\frac{2}{5}+c_3\right)\xi^{3/2}\right)\nonumber\\
&-c_0\partial_{\xi}[\sqrt{\lambda\xi}H_0(\lambda\xi
)]\frac{1}{2}2^{1/4}\lambda (\xi^{5/2}\log\xi
+c_3\xi^{5/2})+O(\lambda^{-\eps}).\nonumber
\end{align}
Evaluating at $\xi=\lambda^{-1/2}$, one obtains that the third and
fourth terms in (\ref{40}) are $O(\lambda^{\frac{1}{2}-\eps})$, and
thus error terms.  Thus,
\begin{align*}
2^{-1/4}c_0^{-1}a_+'(\lambda )=&\left(\frac{1}{2}\lambda^{-1/2}(1+ic_1\log (\lambda\xi )+i\varkappa )+ic_1\lambda^{-1/2}\right)\left(\frac{1}{2}(c_2+\log\xi )+1\right)\\
&-\left(\frac{1}{4}\lambda^{-1/2}(1+ic_1\log (\lambda\xi )+i\varkappa )+ic_1\lambda^{-1/2}\right)(\log\xi +c_2)\\
&+O(\lambda^{-\eps})\\
=&\frac{1}{2}\lambda^{-1/2}(1+ic_1\log (\lambda\xi )+i\varkappa )+ic_1\lambda^{-1/2}\\
&-\frac{ic_1}{2}\lambda^{-1/2}(\log\xi +c_2)+O(\lambda^{-\eps})\\
=&\frac{1}{2}\lambda^{-1/2}(1+ic_1\log\lambda+i\varkappa+2ic_1-ic_1c_2)+O(\lambda^{-\eps})\\
=&\frac{1}{2}\lambda^{-1/2}(1+ic_3+2ic_1+ic_1\log\lambda
)+O(\lambda^{-\eps}).
\end{align*}
Similarly,
\begin{align*}
2^{-1/4}c_0^{-1}b_+'(\lambda )=&\left(\frac{1}{2}\lambda^{-1/2}(1+ic_1\log (\lambda\xi )+i\varkappa )+ic_1\lambda^{-1/2}\right)\left(\frac{1}{2}\right)\\
&-\left(\frac{1}{4}\lambda^{-1/2}(1+ic_1\log (\lambda\xi )+i\varkappa +ic_1\lambda^{-1/2}\right)+O(\lambda^{-\eps})\\
=&-\frac{1}{2}ic_1\lambda^{-1/2}+O(\lambda^{-\eps}),
\end{align*}
as claimed.
\end{proof}

Having this and an explicit expression for the Wronskian
$W(\lambda)$ in terms of $a_ +$ and $b_+$, we obtain

\begin{cor}\label{cor7}
For small $\lambda$,
\begin{align*}
W(\lambda )&=2\lambda \left( 1+ic_3+i\frac{2}{\pi}\log\lambda \right)+O(\lambda^{\frac{3}{2}-\eps})\\
W'(\lambda )&=2\left(
1+ic_3+i\frac{2}{\pi}+i\frac{2}{\pi}\log\lambda \right)
+O(\lambda^{\frac{1}{2}-\eps})
\end{align*}
with $\eps >0$ arbitrary.
\end{cor}

\begin{proof}
Follows immediately from
$$W(\lambda )=-2a_+b_+(\lambda )$$
and (\ref{36}), (\ref{39}).
\end{proof}

To estimate the oscillatory integral \eqref{12} for $\vert
\xi\lambda\vert >1$ we also need the following lemma

\begin{lem}\label{lemma11}
Let $m_+(\xi,\lambda ):=e^{-i\lambda\xi}f_+(\xi,\lambda )$.  Then,
provided $\lambda >0$ is small and $\lambda\xi >1$,
\begin{align}\label{41}
|m_+(\xi,\lambda )-1|&\les (\lambda\xi )^{-1}\\
|\partial_{\lambda}m_+(\xi,\lambda )|&\les
\lambda^{-2}\xi^{-1}\nonumber
\end{align}
\end{lem}

\begin{proof}
From (\ref{24}), and with $m_0(\xi,\lambda
)=e^{-i\lambda\xi}f_0(\xi,\lambda )$, \beeq\label{42}
m_+(\xi,\lambda )=m_0(\xi,\lambda
)+\int_{\xi}^{\infty}\widetilde{G}_0(\xi,\eta;\lambda )V_1(\eta
)m_+(\eta,\lambda )\, d\eta \eneq where \beeq\label{43}
\widetilde{G}_0(\xi,\eta;\lambda )=\frac{m_0(\xi,\lambda
)\overline{m_0(\eta,\lambda )}-e^{-2i(\xi
-\eta)\lambda}\overline{m_0(\xi,\lambda )}m_0(\eta,\lambda
)}{-2i\lambda} \eneq Now, by asymptotic properties of the Hankel
functions,
$$m_0(\xi,\lambda )=1+O((\xi\lambda )^{-1})$$
where the $O$-term behaves like a symbol.\\

Inserting this bound into (\ref{43}) yields
$$|\widetilde{G}_0(\xi,\eta;\lambda )|\les\eta$$
provided $\eta >\xi >\lambda^{-1}.$ Thus, from (\ref{42}),
$$|m_+(\xi,\lambda )-m_0(\xi,\lambda )|\les\xi^{-1}$$
and thus, for all $\xi\lambda >1$,
$$|m_+(\xi,\lambda )-1|\les (\xi\lambda )^{-1}$$
as claimed.\\
Next, one checks that for $\eta >\xi >\lambda^{-1}$,
$$|\partial_{\lambda}\widetilde{G}_0(\xi,\eta;\lambda )|\les\frac{\eta}{\lambda}.$$
Thus, for all $\lambda\xi >1$,
\begin{align*}
|\partial_{\lambda}m_+(\xi,\lambda )|\les &\lambda^{-2}\xi^{-1}+\int_{\xi}^{\infty}|\partial_{\lambda}\widetilde{G}_0(\xi,\eta;\lambda )|\eta^{-3}\, d\eta\\
& + \int_{\xi}^{\infty}\eta^{-2}|\partial_{\lambda}m_+(\eta\lambda )|\, d\eta\\
\les & \lambda^{-2}\xi^{-1}+\lambda^{-1}\xi^{-1}+\int_{\xi}^{\infty}\eta^{-2}|\partial_{\lambda}m_+(\eta,\lambda )|\, d\eta\\
\les & \lambda^{-1}(\lambda\xi )^{-1},
\end{align*}
as claimed.
\end{proof}

We now commence with proving (\ref{12}) for small energies.  Thus,
let $\chi$ be a smooth cut-off function to small energies, i.e.,
$\chi (\lambda )=1$ for small $|\lambda |$ and $\chi$ vanishes
outside a small interval around zero.  In addition, we introduce the
cut-off functions $\chi_{[|\xi\lambda |<1]}$ and $\chi_{[|\xi\lambda
|>1]}$ which form a partition of unity adapted to these intervals.
It will suffice to consider the case $\xi >\xi '$ in (\ref{12}).

\begin{lem}\label{lemma12}
For all $t>0$ \beeq\label{44} \sup_{\xi,\xi
'}\bigg|\int_0^{\infty}e^{it\lambda^2}\lambda\chi (\lambda
)\chi_{[|\xi\lambda |<1,|\xi '\lambda
|<1]}(\langle\xi\rangle\langle\xi
'\rangle)^{-1/2}\mathrm{Im}\left[\frac{f_+(\xi,\lambda )f_-(\xi
',\lambda )}{W(\lambda )}\right]\, d\lambda \bigg|\les t^{-1} \eneq
\end{lem}

\begin{proof}
We first observe the following:
\begin{align*}
\text{Im}&\left[\frac{f_+(\xi,\lambda )f_-(\xi ',\lambda )}{W(\lambda )}\right]\\
&=\text{Im}\left[\frac{(a_+(\lambda )u_0(\xi,\lambda )+b_+(\lambda )u_1(\xi,\lambda ))(a_+(\lambda )u_0(\xi ',\lambda )-b_+(\lambda )u_1(\xi ',\lambda ))}{-2a_+b_+(\lambda )}\right]\\
& = -\frac{1}{2}\text{Im}\left(\frac{a_+}{b_+}(\lambda
)\right)u_0(\xi,\lambda )u_0(\xi ',\lambda
)+\frac{1}{2}\text{Im}\left(\frac{b_+}{a_+}(\lambda
)\right)u_1(\xi,\lambda )u_1(\xi',\lambda).
\end{align*}
Further, by (\ref{36}),
\begin{align*}
-\frac{1}{2}\text{Im}\left(\frac{a_+}{b_+}(\lambda )\right)&=\frac{\pi}{2\sqrt{2}\,}\text{Re}\left[\frac{1+ic_1\log\lambda +ic_3+O(\lambda^{\frac{1}{2}-\eps})}{1+O(\lambda^{\frac{1}{2}-\eps})}\right]\\
&=O(\lambda^{\frac{1}{2}-\eps})+\frac{\pi}{2\sqrt{2}\,}
\end{align*}
and by Corollary \ref{cor6}, the $O$-term can be formally
differentiated, i.e.,
$$\frac{d}{d\lambda}\left\{-\frac{1}{2}\text{Im}\left(\frac{a_+}{b_+}(\lambda )\right)\right\}=O(\lambda^{-\frac{1}{2}-\eps}).$$
Similarly,
$$\frac{1}{2}\text{Im}\left(\frac{b_+}{a_+}(\lambda )\right)=-\frac{\pi}{2\sqrt{2}\,}\,\frac{1}{1+c_1^2\log^2\lambda}+O(\lambda^{\frac{1}{2}-\eps})$$
which can again be formally differentiated.\\
By the estimates of Corollary \ref{cor4}, provided $|\xi\lambda
|+|\xi '\lambda |<1$,
$$|u_0(\xi,\lambda )u_0(\xi ',\lambda )|\les\sqrt{\langle\xi\rangle\langle\xi '\rangle}$$
and
\begin{align*}
|\partial_{\lambda}[u_0(\xi,\lambda )u_0(\xi ',\lambda )]|&\les \lambda \left(\langle\xi\rangle^{5/2}\langle\xi '\rangle^{1/2}+\langle\xi '\rangle^{5/2}\langle\xi\rangle^{1/2}\right)\\
&\les \lambda\sqrt{\langle\xi\rangle\langle\xi
'\rangle}(\langle\xi\rangle^2+\langle\xi '\rangle^2).
\end{align*}
Similarly,
$$|u_1(\xi,\lambda )u_1(\xi ',\lambda )|\les\sqrt{\langle\xi\rangle\langle\xi '\rangle}(\log\lambda )^2$$
and
$$|\partial_{\lambda}[u_1(\xi,\lambda )u_1(\xi ',\lambda )]|\les\lambda\sqrt{\langle\xi\rangle\langle\xi '\rangle}(\langle\xi\rangle^2+\langle\xi '\rangle^2)(\log\lambda )^2.$$
Hence, integrating by parts in (\ref{44}) yields
\begin{align*}
(\ref{44})&\les  \\
t^{-1}&\int_0^{\infty}\bigg|\partial_{\lambda} [\chi (\lambda )\chi_{[|\xi\lambda |,|\xi '\lambda |<1]}
(\langle\xi\rangle\langle\xi '\rangle )^{-1/2}\text{Im}\left(\frac{a_+}{b_+}(\lambda )\right)u_0(\xi,\lambda )
u_0(\xi ',\lambda )]\bigg|\, d\lambda \\
+ &t^{-1}\int_0^{\infty}\bigg|\partial_{\lambda} [\chi (\lambda )\chi_{[|\xi\lambda |,|\xi '\lambda |<1]}
(\langle\xi\rangle\langle\xi '\rangle )^{-1/2}\text{Im}\left(\frac{b_+}{a_+}(\lambda )\right)
u_1(\xi,\lambda )u_1(\xi ',\lambda )]\bigg|\, d\lambda \\
\les&\, 1+\int_0^{\infty}\chi (\lambda )\lambda(\langle\xi\rangle^2+\langle\xi '\rangle^2)
\chi_{[|\xi\lambda |,|\xi '\lambda |<1]}\, d\lambda\\
\les &\, 1,
\end{align*}
and the lemma is proved.
\end{proof}
Next, we consider the case $|\xi\lambda |>1$ and $|\xi '\lambda
|>1$.  With the convention that $f_{\pm}(\xi,-\lambda
)=\overline{f_{\pm}(\xi,\lambda )}$ we can remove the imaginary part
in (\ref{12}) and integrate $\lambda$ over the whole axis.  We shall
follow this convention hence forth. To estimate the oscillatory
integrals, we shall repeatedly use the following version of
stationary
phase, see Lemma 2 in \cite{Sch} for the proof.\\
\begin{lem}\label{lemma12'}
Let $\phi (0)=\phi '(0)=0$ and $1\leq\phi ''\leq C$.  Then
\beeq\label{46} \bigg|\int_{-\infty}^{\infty}e^{it\phi (x)}a(x)\,
dx\bigg|\les\delta^2\left\{\int\frac{|a(x)|}{\delta^2+|x|^2}\,
dx+\int_{|x|>\delta}\frac{|a'(x)|}{|x|}\, dx\right\} \eneq
where $\delta =t^{-1/2}$.\\
\end{lem}

\noindent Using Lemma \ref{lemma12'} we can prove the following:
\begin{lem}\label{lemma13}
For all $t>0$ \beeq\label{45} \sup_{\xi >0>\xi
'}\bigg|\int_{-\infty}^{\infty}e^{it\lambda^2}\lambda\chi (\lambda
)\chi_{[|\xi\lambda |>1,|\xi '\lambda
|>1]}(\langle\xi\rangle\langle\xi
'\rangle)^{-1/2}\frac{f_+(\xi,\lambda )f_-(\xi ',\lambda
)}{W(\lambda )}\, d\lambda \bigg|\les t^{-1} \eneq
\end{lem}

\begin{proof}
Writing $f_+(\xi,\lambda )=e^{i\xi\lambda}m_+(\xi,\lambda )$,
$f_-(\xi,\lambda )=e^{-i\xi\lambda}m_-(\xi,\lambda )$ as in Lemma
\ref{lemma11}, we express (\ref{45}) in the form \beeq\label{47}
\bigg|\int_{-\infty}^{\infty}e^{it\phi (\lambda )}a(\lambda )\,
d\lambda \bigg|\les t^{-1} \eneq where $\xi > 0>\xi '$ are fixed,
$\phi (\lambda ):=\lambda^2+\frac{\lambda}{t}(\xi -\xi ')$, and
$$a(\lambda )=\lambda\chi (\lambda )\chi_{[|\xi\lambda |>1,|\xi '\lambda |>1]}(\langle\xi\rangle\langle\xi '\rangle )^{-1/2}\frac{m_+(\xi,\lambda )m_-(\xi ',\lambda )}{W(\lambda )}.$$
Let $\lambda_0=-\frac{\xi -\xi '}{2t}$. We have the bounds
\beeq\label{48} |a(\lambda )|\les(\langle\xi\rangle\langle\xi
'\rangle )^{\frac{-1}{2}}\chi (\lambda )\chi_{[|\xi\lambda |>1,|\xi
'\lambda |>1]}. \eneq By Corollary \ref{cor7}, for small $|\lambda
|$
$$\bigg|\left(\frac{\lambda}{W(\lambda )}\right)'\bigg|\les\frac{1}{|\lambda |(\log |\lambda |)^2}$$
and by Lemma \ref{lemma11}, for $|\xi\lambda |>1$, $|\xi '\lambda
|>1$,
$$|\partial_{\lambda}[m_+(\xi,\lambda )m_-(\xi ',\lambda )]|\les\lambda^{-2}(\xi^{-1}+|\xi '|^{-1}).$$
Hence, \beeq\label{49} |a'(\lambda
)|\les(\langle\xi\rangle\langle\xi '\rangle )^{-1/2}\chi (\lambda
)\chi_{[|\xi\lambda |>1,|\xi '\lambda |>1]}\left\{\frac{|\lambda
|^{-1}}{|\log\lambda |^2}+\lambda^{-2}(\xi^{-1}+|\xi
'|^{-1})\right\}. \eneq
We will need to consider three cases in order to prove (\ref{47}) via (\ref{46}), depending on where $\lambda_0$ falls relative to the support of $a$.\\

\underline{Case 1:} $|\lambda_0 |\les 1$, $|\lambda_0|\gtr |\xi |^{-1}+|\xi '|^{-1}$.\\

Note that the second inequality here implies that
$$\frac{\xi +|\xi '|}{t}\gtr\frac{\xi +|\xi '|}{\xi |\xi '|}$$
or
$$1\gtr\frac{t}{\xi |\xi '|}.$$
Furthermore, we remark that $a\equiv 0$ unless $\xi\gtr 1$ and $|\xi '|\gtr 1$.\\
Starting with the first integral on the right-hand side of
(\ref{46}) we conclude from (\ref{48}) that
$$\int\frac{|a(\lambda )|}{|\lambda -\lambda_0|^2+\delta^2}\les(\langle\xi\rangle\langle\xi '\rangle )^{-1/2}t^{1/2}\les 1.$$
From the second integral we obtain from (\ref{49}) that
\begin{align*}
\int_{|\lambda -\lambda_0|>\delta}\frac{|a'(\lambda )|}{|\lambda -\lambda_0|}\, d\lambda \les &(\langle\xi\rangle\langle\xi '\rangle )^{-1/2}\delta^{-1}\int\frac{\chi (\lambda )\, d\lambda}{|\lambda |(\log |\lambda |)^2}\\
&+(\langle\xi\rangle\langle\xi '\rangle
)^{-1/2}(\langle\xi\rangle^{-1}
+\langle\xi '\rangle^{-1})\delta^{-1}\int_{\lambda >\xi^{-1}+|\xi '|^{-1}}\frac{d\lambda}{\lambda^2}\\
&\les\sqrt{\frac{t}{\langle\xi\rangle\langle\xi '\rangle}}\\
&\les 1.
\end{align*}
\\

\underline{Case 2:} $|\lambda_0 |\les 1$, $|\lambda_0|<< \langle\xi\rangle^{-1}+\langle\xi '\rangle^{-1}$.\\

Then $|\lambda -\lambda_0|\sim |\lambda |$ on the support of $a$,
which implies that
\begin{align*}
\int\frac{|a(\lambda )|}{|\lambda -\lambda_0|^2+t^{-1}}\, d\lambda &\les (\langle\xi\rangle\langle\xi '\rangle )^{-1/2}\int_{\lambda >\xi^{-1}+|\xi '|^{-1}}\frac{d\lambda}{\lambda^2}\\
&\les\frac{\sqrt{\xi |\xi '|}}{\xi +|\xi '|}\\
&\les 1,
\end{align*}
and also
\begin{align*}
\int_{|\lambda -\lambda_0|>\delta}\frac{|a'(\lambda )|}{|\lambda -\lambda_0|}\, d\lambda \les&(\langle\xi\rangle\langle\xi '\rangle )^{-1/2}(\int_{\lambda >\xi^{-1}+|\xi '|^{-1}}\frac{d\lambda}{\lambda^2(\log |\lambda |)^2}\\
&+\int_{\lambda >\xi^{-1}+|\xi '|^{-1}}\frac{d\lambda}{\lambda^3}(\xi^{-1}+|\xi '|^{-1}))\\
\les&\frac{\sqrt{\xi |\xi '|}}{\xi +|\xi '|}\\
&\les 1.
\end{align*}
\\

\underline{Case 3:} $|\lambda_0 |>> 1$, $|\lambda_0|\gtr\xi^{-1}+|\xi ' |^{-1}$.\\

In this case, $|\lambda - \lambda_0 |\sim |\lambda_0|>>1$. Thus,
$$\int\frac{|a(\lambda )|}{|\lambda -\lambda_0|^2+t^{-1}}\, d\lambda \les (\langle\xi\rangle\langle\xi '\rangle )^{-1/2}\frac{1}{\lambda_0^2+t^{-1}}\les 1$$
as well as, see (\ref{49}),
\begin{align*}
\int_{|\lambda -\lambda_0|>\delta}\frac{|a'(\lambda )|}{|\lambda -\lambda_0|}\, d\lambda \les& \int\frac{\chi (\lambda )}{|\lambda |(\log |\lambda |)^2}d\lambda \frac{(\langle\xi\rangle\langle\xi '\rangle )^{-1/2}}{\lambda_0}\\
& +\int\frac{1}{\lambda^2}\chi_{[|\lambda |>\xi^{-1}+|\xi '|^{-1}]}\frac{d\lambda}{\lambda_0}\:\frac{\xi +|\xi '|}{(\xi |\xi '|)^{3/2}}\\
&\les 1,
\end{align*}
and the lemma is proved.
\end{proof}

Now we turn to the estimate of the oscillatory integral for the case
$\xi\lambda
>1$ and $|\xi '\lambda |<1$.

\begin{lem}\label{lemma14}
For all $t>0$ \beeq\label{50} \sup_{\xi >0>\xi
'}\bigg|(\langle\xi\rangle\langle\xi '\rangle
)^{-1/2}\int_{-\infty}^{\infty}e^{it\lambda^2}\frac{\lambda\chi
(\lambda )}{W(\lambda )}\chi_{[\xi\lambda >1,|\xi '\lambda
|<1]}f_+(\xi,\lambda )f_-(\xi ',\lambda )\, d\lambda \bigg|\les
t^{-1} \eneq and similarly with $\chi_{[|\xi\lambda |<1,\xi '\lambda
<-1]}$.
\end{lem}

\begin{proof}
As before, we write $f_+(\xi,\lambda
)=e^{i\xi\lambda}m_+(\xi,\lambda )$.  But because of $|\xi '\lambda
|<1$ we use the representation
$$f_-(\xi ',\lambda )=a_-(\lambda )u_0(\xi,\lambda )+b_-(\lambda )u_1(\xi,\lambda ).$$
In particular,
$$|f_-(\xi ',\lambda )|\les\sqrt{|\lambda |\langle\xi '\rangle}\big|\log |\lambda |\big|.$$
Moreover, from (\ref{34}) and (\ref{39}),
$$|\partial_{\lambda}f_-(\xi ',\lambda )|\les\langle\xi '\rangle^{1/2}|\lambda |^{-1/2}\big|\log |\lambda |\big|$$
provided $|\xi '\lambda |<1$.\\
We apply (\ref{46}) with $\phi (\lambda
)=\lambda^2+\frac{\xi}{t}\lambda$ and
$$a(\lambda )=\frac{\lambda\chi (\lambda )}{W(\lambda )}(\langle\xi\rangle\langle\xi '\rangle )^{-1/2}\chi_{[\xi\lambda >1,|\xi '\lambda |<1]}m_+(\xi,\lambda )f_-(\xi ',\lambda ).$$
By the preceding, \beeq\label{51} |a(\lambda )|\les\frac{|\lambda
|^{1/2}}{\sqrt{\langle\xi\rangle}}\chi (\lambda )\chi_{[\xi\lambda
>1,|\xi '\lambda |<1]} \eneq and \beeq\label{52} |a'(\lambda )|\les
(|\lambda |\langle\xi\rangle )^{-1/2}\chi (\lambda
)\chi_{[\xi\lambda >1,|\xi '\lambda |<1]}.\eneq

\underline{Case 1:} $|\lambda_0|\les 1$, $|\xi\lambda_0|\gtr 1$.\\

Note in particular $|\xi |\gtr 1$. Here $\lambda_0=-\frac{\xi}{2t}$.
By (\ref{51}),
\begin{align*}
\int\frac{|a(\lambda )|}{|\lambda -\lambda_0|^2+t^{-1}}d\lambda &\les\langle\xi\rangle^{-1/2}\int\frac{\sqrt{|\lambda |}}{|\lambda -\lambda_0|^2+t^{-1}}d\lambda\\
&\les\langle\xi\rangle^{-1/2}|\lambda_0|^{1/2}\int\frac{d\lambda}{|\lambda -\lambda_0|^2+t^{-1}}+\langle\xi\rangle^{-1/2}\int\frac{|\lambda |^{1/2}}{|\lambda |^2+t^{-1}}d\lambda\\
&\les\langle\xi\rangle^{-1/2}t^{1/2}\left(\frac{\xi}{t}\right)^{1/2}+\langle\xi\rangle^{-1/2}t^{1/4}\\
&\les 1
\end{align*}
Here we used that $|\xi\lambda_0|=\frac{\xi^2}{2t}\gtr 1$.\\
Next, write via (\ref{52}) \beeq\label{53} \int_{|\lambda
-\lambda_0|>\delta}\frac{|a'(\lambda )|}{|\lambda
-\lambda_0|}d\lambda\les\langle\xi\rangle^{-\frac{1}{2}}\int_{|\lambda
-\lambda_0|>\delta}\frac{1}{|\lambda |^{\frac{1}{2}}|\lambda
-\lambda_0|}\chi_{[\xi\lambda >1,|\xi '\lambda |<1]}\, d\lambda.
\eneq Distinguish the cases $\frac{1}{10}|\lambda |>|\lambda
-\lambda_0|$ and $\frac{1}{10}|\lambda |\leq |\lambda -\lambda_0|$
in the integral on the right-hand side.  This yields
\begin{align*}
(\ref{53})\les&\langle\xi\rangle^{-1/2/}\int_{|\lambda -\lambda_0|>\delta}\frac{d\lambda}{|\lambda -\lambda_0|^{3/2}}+\langle\xi\rangle^{-1/2}\int_{|\lambda |\les |\lambda_0|}\frac{d\lambda}{|\lambda |^{1/2}}|\lambda_0|^{-1}\\
&+\langle\xi\rangle^{-1/2}\int_{|\lambda |>|\lambda_0|}\frac{d\lambda}{|\lambda |^{3/2}}\\
\les &\langle\xi\rangle^{-1/2}\delta^{-1/2}+\langle\xi\rangle^{-1/2}|\lambda_0|^{-1/2}\\
\les &\left(\frac{t}{\xi^2}\right)^{1/4}+|\xi\lambda_0|^{-1/2}\\
\les & 1.
\end{align*}
\\

\underline{Case 2:} $|\lambda_0|\les 1$, $|\xi\lambda_0 |<<1$.\\

In that case, $|\lambda -\lambda_0|\sim|\lambda |$ on the support of
$a$.  Consequently,
$$\int\frac{|a(\lambda )|}{|\lambda -\lambda_0|^2+t^{-1}}d\lambda\les\langle\xi\rangle^{-\frac{1}{2}}\int_{|\xi |^{-1}}^{\infty}|\lambda |^{-\frac{3}{2}}\, d\lambda\les 1.$$
Moreover,
$$\int_{|\lambda -\lambda_0|>\delta}\frac{|a'(\lambda )|}{|\lambda -\lambda_0|}d\lambda\les\int_{|\xi |^{-1}}^{\infty}\frac{(|\lambda |\langle\xi\rangle )^{-\frac{1}{2}}}{|\lambda |}d\lambda \les 1.$$
\\

\underline{Case 3:} $|\lambda_0|>>1$.\\

In that case, $|\lambda -\lambda_0|\sim |\lambda_0|$ on $\supp (a)$.
Since $|a(\lambda )|\les 1$ by (\ref{51}), it follows that
$$\int\frac{|a(\lambda )|}{|\lambda -\lambda_0|^2+t^{-1}}d\lambda\les 1.$$
Similarly, since $|a'(\lambda )|\les (\xi |\lambda
|)^{-\frac{1}{2}}$, it follows that
$$\int_{|\lambda -\lambda_0|>\delta}\frac{|a'(\lambda )|}{|\lambda -\lambda_0|}d\lambda\les\int\frac{(|\lambda |\langle\xi\rangle )^{-\frac{1}{2}}}{|\lambda_0|}\chi (\lambda )\, d\lambda\les 1.$$
This proves (\ref{50}).  The other case $\chi_{[|\xi\lambda |<1,\xi
'\lambda <-1]}$ is treated in an analogous fashion.
\end{proof}

The remaining cases for the small energy part of (\ref{12}) are $\xi
>\xi '>|\lambda |^{-1}$ and $\xi '<\xi <-|\lambda |^{-1}$.  By
symmetry it will suffice to treat the former case.  As usual, we
need to consider reflection and transmission coefficients, therefore
we write \beeq\label{54} f_-(\xi,\lambda )=\alpha_-(\lambda
)f_+(\xi,\lambda )+\beta_-(\lambda )\overline{f_+(\xi,\lambda )}.
\eneq Then, with $W(\lambda )=W(f_+(\cdot,\lambda
),f_-(\cdot,\lambda ))$,
$$W(\lambda )=\beta_-(\lambda )W(f_+(\cdot,\lambda ),\overline{f_+(\cdot,\lambda )})=-2i\lambda\beta_-(\lambda )$$
and
\begin{align*}
W(f_-(\cdot,\lambda ),\overline{f_+(\cdot,\lambda )})&=\alpha_-(\lambda )W(f_+(\cdot,\lambda )\overline{f_+(\cdot,\lambda )})\\
&=-2i\lambda\alpha_-(\lambda ).
\end{align*}
Thus, when $\lambda >0$ is small, \beeq\label{55} \beta_-(\lambda
)=i\left(1+ic_3+i\frac{2}{\pi}\log\lambda \right)+O(|\lambda
|^{\frac{1}{2}-\eps}) \eneq and
\begin{align}\label{56}
\alpha_-(\lambda )&=\frac{1}{-2i\lambda}W\left(a_+(\lambda )u_0(\cdot,\lambda )-b_+(\lambda )u_1(\cdot,\lambda ),\overline{a_+}(\lambda )u_0(\cdot,\lambda )+\overline{b_+}(\lambda )u_1(\cdot,\lambda )\right)\nonumber\\
&=\frac{1}{-2i\lambda}\left(a_+\overline{b_+}(\lambda )+\overline{a_+}(\lambda )b_+(\lambda )\right)\nonumber\\
&=\frac{i}{\lambda}\text{Re}(a_+\overline{b_+}(\lambda))\nonumber\\
&=\frac{i}{\lambda}\text{Re}\left(-i|c_0|^2c_1\lambda(1+ic_1\log\lambda +ic_3)+O(\lambda^{\frac{3}{2}-\eps})\right)\nonumber\\
&=i\left(\frac{2}{\pi}\log\lambda
+c_3\right)+O(\lambda^{\frac{1}{2}-\eps}).
\end{align}
In passing, we remark that $1+|\alpha_-|^2=|\beta_-|^2$.  Finally,
it follows from Corollary \ref{cor6} that the $O$-terms can be
differentiated once in $\lambda$; they then become
$O(\lambda^{-\frac{1}{2}-\eps})$, $\eps >0$ arbitrary.

\begin{lem}\label{lemma15}
For any $t>0$ \beeq\label{57} \sup_{\xi >\xi
'>0}\bigg|(\langle\xi\rangle\langle\xi '\rangle )^{-\frac{1}{2}}\int
e^{it\lambda^2}\frac{\lambda\chi (\lambda )}{W(\lambda )}\chi_{[|\xi
'\lambda |>1]}f_+(\xi,\lambda )f_-(\xi ',\lambda )\, d\lambda
\bigg|\les t^{-1} \eneq and similarly for $\sup_{\xi ' <\xi <0}$ and
$\chi_{[|\xi\lambda |>1]}$.
\end{lem}

\begin{proof}
Using (\ref{54}), we reduce (\ref{57}) to two estimates:
\begin{align}\label{58}
\sup_{\xi >\xi '>0}&\bigg|(\langle\xi\rangle\langle\xi '\rangle )^{-\frac{1}{2}}\int e^{it\lambda^2}\frac{\lambda\chi (\lambda )}{W(\lambda )}\chi_{[\xi '|\lambda |>1]}e^{i\lambda (\xi +\xi ')}m_+(\xi,\lambda )m_+(\xi ',\lambda )\alpha_-(\lambda )\, d\lambda \bigg|\\
&\les t^{-1}\nonumber
\end{align}
and
\begin{align}\label{59}
\sup_{\xi >\xi '>0}&\bigg|(\langle\xi\rangle\langle\xi '\rangle )^{-\frac{1}{2}}\int e^{it\lambda^2}e^{i\lambda (\xi -\xi ')}\frac{\lambda\chi (\lambda )}{W(\lambda )}\beta_-(\lambda )\chi_{[\xi '|\lambda |>1]}m_+(\xi,\lambda )\overline{m_+(\xi ',\lambda )}\, d\lambda \bigg|\\
&\les t^{-1}\nonumber
\end{align}
We apply (\ref{46}) to (\ref{58}) with fixed $\xi >\xi '>0$ and
\begin{align*}
\phi (\lambda )&=\lambda^2+\frac{\lambda}{t}(\xi +\xi '),\\
a(\lambda )&=(\langle\xi\rangle\langle\xi '\rangle
)^{-\frac{1}{2}}\frac{\lambda\chi (\lambda )}{W(\lambda )}\chi_{[\xi
'|\lambda |>1]}\alpha_-(\lambda )m_+(\xi,\lambda )m_+(\xi ',\lambda
).
\end{align*}
Then from (\ref{56}), \beeq\label{60} |a(\lambda
)|\les(\langle\xi\rangle\langle\xi '\rangle )^{-\frac{1}{2}}\chi
(\lambda )\chi_{[\xi '|\lambda |>1]} \eneq and from our derivative
bounds on $W$, $\alpha_-$, and $m_+(\xi,\lambda)$, see (\ref{41})
for the latter, we conclude that \beeq\label{61} |a'(\lambda )|\les
|\lambda |^{-1}(\langle\xi\rangle\langle\xi '\rangle
)^{-\frac{1}{2}}\chi (\lambda )\chi_{[\xi '|\lambda |>1]}. \eneq
\\

\underline{Case 1:} Suppose $|\lambda_0|\les 1$ and $|\xi '\lambda_0 |>1$, where $\lambda_0=-\frac{\xi +\xi '}{2t}$.  Note $\xi >\xi '\gtr 1$.\\
Then
\begin{align*}
\int\frac{|a(\lambda )|}{|\lambda -\lambda_0|^2+t^{-1}}d\lambda &\les (\langle\xi\rangle\langle\xi '\rangle )^{-\frac{1}{2}}\int\frac{d\lambda}{|\lambda -\lambda_0|^2+t^{-1}}\\
&\les\sqrt{\frac{t}{\xi\xi '}}\les 1
\end{align*}
since $|\xi '\lambda_0|\sim\frac{\xi\xi '}{t}>1$.\\
As for the derivative term in (\ref{46}), we infer from (\ref{61})
that \beeq\label{62} \int_{|\lambda
-\lambda_0|>\delta}\frac{|a'(\lambda )|}{|\lambda -\lambda_0
|}d\lambda\les(\langle\xi\rangle\langle\xi '\rangle
)^{-\frac{1}{2}}\int_{|\lambda
-\lambda_0|>\delta}\frac{d\lambda}{|\lambda ||\lambda -\lambda_0
|}\chi_{[|\lambda\xi '|>1]} \eneq Again, we need to distinguish
between $|\lambda -\lambda_0|>\frac{1}{10}|\lambda_0|$ and $|\lambda
-\lambda_0|<\frac{1}{10}|\lambda_0|$.  Thus, since $\xi\xi '>t$,
\begin{align*}
(\ref{62})&\les (\langle\xi\rangle\langle\xi '\rangle )^{-\frac{1}{2}}\int_{1/\xi '}^{\infty}\frac{d\lambda}{\lambda^2}+(\langle\xi\rangle\langle\xi '\rangle )^{-1/2}|\lambda_0|^{-1}\log \left(t^{1/2}|\lambda_0|\right)\\
&\les 1+\frac{t^{\frac{1}{2}}}{\xi}\log \left(\frac{\xi}{t^{1/2}}\right)\\
&\les 1
\end{align*}
since also $\xi^2>t$.\\

\underline{Case 2:} $|\lambda_0|\les 1$, $|\lambda_0|<<\frac{1}{\xi '}$.\\

Then $|\lambda -\lambda_0|\sim|\lambda |$ on the support of
$a(\lambda )$.  Hence,
$$\int\frac{|a(\lambda )|}{|\lambda -\lambda_0|^2+t^{-1}}d\lambda\les (\langle\xi\rangle\langle\xi '\rangle )^{-\frac{1}{2}}\int_{1/\xi '}^{\infty}\frac{d\lambda}{\lambda^2}\les\sqrt{\frac{\xi '}{\langle\xi\rangle}}<1$$
and
$$\int_{|\lambda -\lambda_0|>\delta}\frac{|a'(\lambda )|}{|\lambda -\lambda_0|}d\lambda\les (\langle\xi\rangle\langle\xi '\rangle )^{-\frac{1}{2}}\int_{1/\xi '}^{\infty}\frac{d\lambda}{\lambda^2}<1.$$
\\

\underline{Case 3:} $|\lambda_0|>>1$, $|\lambda_0|\gtr\frac{1}{\xi '}$.\\

Then $|\lambda -\lambda_0|\sim |\lambda_0|$ on $\supp (a)$.
Therefore, $|a (\lambda )|\les 1$ implies that
$$\int\frac{|a(\lambda )|}{|\lambda -\lambda_0|^2+t^{-1}}d\lambda\les 1$$
and
\begin{align*}
\int_{|\lambda -\lambda_0|>\delta}\frac{|a'(\lambda )|}{|\lambda -\lambda_0|}d\lambda&\les (\langle\xi\rangle\langle\xi '\rangle )^{-\frac{1}{2}}|\lambda_0|^{-1}\int_{\frac{1}{\langle\xi '\rangle}}^{1}\frac{d\lambda}{|\lambda |}\\
&\les\frac{1}{|\lambda_0|}\frac{1}{\langle\xi '\rangle}\log\langle\xi '\rangle\\
&\les 1.
\end{align*}
This concludes the proof of (\ref{58}). (\ref{59}) is completely
analogous, as is the case of $\xi '<\xi <0$, $|\xi\lambda |>1$.
\end{proof}

We are done with the contributions of small $\lambda$ in (\ref{12}).
To conclude the proof of Theorem \ref{thm1} it suffices to prove the
following statement.

\begin{lem}\label{lemma16}
For all $t>0$, \beeq\label{63} \sup_{\xi >\xi
'}\bigg|(\langle\xi\rangle\langle\xi '\rangle
)^{-\frac{1}{2}}\int_{-\infty}^{\infty}e^{it\lambda^2}\frac{\lambda
(1-\chi )(\lambda )}{W(\lambda )}f_+(\xi,\lambda )f_-(\xi ',\lambda
)\, d\lambda \bigg|\les t^{-1}. \eneq
\end{lem}

\begin{proof}
We observed above, see (\ref{54}), that $W(\lambda )=-2i\lambda\beta_-(\lambda )$.  Since $|\beta_-(\lambda )|\geq 1$, this implies that $|W|(\lambda )\geq 2|\lambda |$.  In particular, $W(\lambda )\neq 0$ for every $\lambda\neq 0$.\\
In order to prove (\ref{63}), we will need to distinguish the cases $\xi >0>\xi '$, $\xi >\xi '>0$, $0>\xi >\xi '$.  By symmetry, it will suffice to consider the first two.\\

\underline{Case 1:} $\xi >0>\xi '$.\\

In this case we need to prove that \beeq\label{64} \sup_{\xi >0>\xi
'}\bigg|(\langle\xi\rangle\langle\xi '\rangle )^{-\frac{1}{2}}\int
e^{it[\lambda^2+\frac{\xi -\xi '}{t}\lambda]}\frac{\lambda (1-\chi
)(\lambda )}{W(\lambda )}m_+(\xi,\lambda )m_-(\xi ',\lambda )\,
d\lambda \bigg|\les t^{-1}. \eneq Apply (\ref{46}) with $\phi
(\lambda )=\lambda^2+\frac{\xi -\xi '}{t}\lambda$ and
$$a(\lambda )= (\langle\xi\rangle\langle\xi '\rangle )^{-\frac{1}{2}}\frac{\lambda (1-\chi )(\lambda )}{W(\lambda )}m_+(\xi,\lambda )m_-(\xi ',\lambda ).$$
Hence, with $\lambda_0=-\frac{\xi -\xi '}{2t}$,
\begin{align*}
(\ref{64})&\les t^{-1}\left(\int\frac{|a(\lambda )|}{|\lambda -\lambda_0|^2+t^{-1}}d\lambda +\int_{|\lambda -\lambda_0|>\delta}\frac{|a'(\lambda )|}{|\lambda -\lambda_0|}d\lambda\right)\\
&=:t^{-1}(A+B).
\end{align*}
If $|\lambda_0|<<1$, then
$$A\les ||a||_{\infty}\les 1.$$
On the other hand, if $|\lambda_0|\gtr 1$, then $\xi +|\xi '|\gtr t$
so that
$$A\les t^{\frac{1}{2}}||a||_{\infty}\les t^{\frac{1}{2}}(\langle\xi\rangle\langle\xi '\rangle )^{-\frac{1}{2}}\les\sqrt{\frac{t}{\langle\xi\rangle\langle\xi '\rangle}}\les 1.$$
Here we used that
$$\sup_{\xi}\sup_{|\lambda |\gtr 1}|m_{\pm}(\xi,\lambda )|\les 1$$
which follows from the fact that \beeq\label{65} m_+(\xi,\lambda
)=1+\int_{\xi}^{\infty}\frac{1-e^{-2i(\tilde{\xi} -\xi
)\lambda}}{2i\lambda}V(\tilde{\xi} )m_+(\tilde{\xi},\lambda
)d\tilde{\xi} \eneq with
$$V(\tilde{\xi} )=\left(\frac{1}{4}\omega^2+\frac{1}{2}\dot{\omega}\right)(\tilde{\xi})=O(\langle\tilde{\xi}\rangle^{-2}).$$
Moreover, from our assumptions on $r(x)$ we recall that
$$\bigg|\frac{d^l}{d\xi^l}V(\xi )\bigg|\les\langle\xi\rangle^{-2-l},\qquad l\geq 0.$$
We shall need these bounds to estimate $B$ above. From (\ref{65}),
for $\xi\geq 0$
$$m_+(\xi,\lambda )=1+O(\lambda^{-1}\langle\xi\rangle^{-1})$$
as well as for $\xi\geq 0$
\begin{align}\label{66}
\partial^j_{\xi}m_+(\xi,\lambda )&=O(\lambda^{-1}\langle\xi\rangle^{-1}),\qquad j=1,2
\end{align}
\begin{align}\label{67}
\partial_{\lambda}m_+(\xi,\lambda )&=O(\lambda^{-2}\langle\xi\rangle^{-1})
\end{align}
\begin{align}\label{68}
\partial_{\lambda}\partial_{\xi}m_+(\xi,\lambda )&=O(\lambda^{-2}\langle\xi\rangle^{-1})
\end{align}
To verify (\ref{66}), one checks that
\begin{align}\label{69}
\partial_{\xi}m_+(\xi,\lambda )=&\frac{1}{2i\lambda}\int_{\xi}^{\infty}e^{2i(\xi -\tilde{\xi})\lambda}\partial_{\tilde{\xi}}\left[(\xi -\tilde{\xi})V(\tilde{\xi})\right]m_+(\tilde{\xi},\lambda )\, d\tilde{\xi}\\
&+\frac{1}{2i\lambda}\int_{\xi}^{\infty}e^{2i(\xi
-\tilde{\xi})\lambda}(\xi
-\tilde{\xi})V(\tilde{\xi})\partial_{\tilde{\xi}}m_+(\tilde{\xi},\lambda
)\, d\tilde{\xi}.\nonumber
\end{align}
By our estimates on $V$, the integral in (\ref{69}) is
$O(\lambda^{-1}\langle\xi\rangle^{-1})$ and (\ref{66}) follows.  For
(\ref{67}) we compute
\begin{align*}
\partial_{\lambda}m_+(\xi,\lambda )=-&\int_{\xi}^{\infty}\frac{1-e^{2i(\xi -\tilde{\xi})\lambda}}{2i\lambda^2}V(\tilde{\xi})m_+(\tilde{\xi},\lambda )\, d\tilde{\xi}\\
&+\frac{1}{2i\lambda^2}\int_{\xi}^{\infty}e^{2i(\xi -\tilde{\xi})\lambda}\partial_{\tilde{\xi}}\left[(\xi -\tilde{\xi})V(\tilde{\xi})m_+(\tilde{\xi},\lambda )\right]\, d\tilde{\xi}\\
&+\int_{\xi}^{\infty}\frac{1-e^{2i(\xi
-\tilde{\xi})\lambda}}{2i\lambda}V(\tilde{\xi}
)\partial_{\lambda}m_+(\tilde{\xi},\lambda )\, d\tilde{\xi}
\end{align*}
so that
$$\partial_{\lambda}m_+(\xi,\lambda )=O(\lambda^{-2}\langle\xi\rangle^{-1})$$
as claimed.\\

Finally, compute
\begin{align*}
\partial^2_{\xi\lambda}m_+(\xi,\lambda )=\frac{1}{\lambda}&\int_{\xi}^{\infty}e^{2i(\xi -\tilde{\xi})\lambda}V(\tilde{\xi})m_+(\tilde{\xi},\lambda )\, d\tilde{\xi}\\
&+\frac{1}{2i\lambda^2} V(\xi )m_+(\xi,\lambda )+ \frac{1}{2i\lambda}\int_{\xi}^{\infty}e^{2i(\xi -\tilde{\xi})\lambda}\partial_{\tilde{\xi}}[(\xi -\tilde{\xi})Vm_+]\, d\tilde{\xi}\\
&+\frac{1}{2i\lambda^2}\int_{\xi}^{\infty}e^{2i(\xi -\tilde{\xi})\lambda}\partial_{\tilde{\xi}}[V(\tilde{\xi})m_+(\tilde{\xi},\lambda )]\, d\tilde{\xi}\\
&-\int_{\xi}^{\infty}e^{2i(\xi
-\tilde{\xi})\lambda}V(\tilde{\xi})\partial_{\lambda}m_+(\tilde{\xi},\lambda
)\, d\tilde{\xi}
\end{align*}
Integrating by parts in the first, third, and last terms yields the
desired estimate.  As a corollary, we obtain (take $\xi =0$)
\begin{align*}
W(\lambda )&=W(f_+(\cdot,\lambda ),f_-(\cdot,\lambda ))\\
&=m_+(\xi,\lambda )[m'_-(\xi,\lambda )-i\lambda m_-(\xi,\lambda )]-m_-(\xi,\lambda )[m'_+(\xi,\lambda )+i\lambda m_+(\xi,\lambda )]\\
&=-2i\lambda (1+O(\lambda^{-1}))+O(\lambda^{-1})\\
&=-2i\lambda+O(1)
\end{align*}
and
$$W'(\lambda )=-2i+O(\lambda^{-1})$$
as $|\lambda |\to\infty$.\\

Next, we estimate $B$. First, we conclude from our bounds on
$W(\lambda )$ and $m_+(\xi,\lambda )$ as well as $m_-(\xi ',\lambda
)$ that
$$|a'(\lambda )|\les (\langle\xi\rangle\langle\xi '\rangle )^{-\frac{1}{2}}\chi_{[|\lambda |\gtr 1]}|\lambda |^{-1}.$$
Let us first consider the case where $|\lambda_0|\gtr 1$.  Then
\begin{align*}
B&\les(\langle\xi\rangle\langle\xi '\rangle )^{-1/2}\int_{\binom{|\lambda -\lambda_0|>\delta}{|\lambda |\gtr 1}}\frac{d\lambda}{|\lambda ||\lambda -\lambda_0|}\\
&\les (\langle\xi\rangle\langle\xi '\rangle )^{-1/2}\left\{\int_1^{\infty}\frac{d\lambda}{\lambda^2}+\frac{1}{|\lambda_0|}\int_{\frac{|\lambda_0|}{5}>|\lambda -\lambda_0 |>\delta}\frac{d\lambda}{|\lambda -\lambda_0|}\right\}\\
&\les 1+\sqrt{\frac{t}{\langle\xi\rangle\langle\xi '\rangle}}\frac{1}{|\lambda_0|t^{1/2}}\log_+(\lambda_0t^{1/2})\\
&\les 1
\end{align*}
Here we used that $\frac{t}{\langle\xi\rangle\langle\xi '\rangle}\les 1$ which follows from $|\lambda_0|\gtr 1$.\\
If $|\lambda_0|<<1$, then $|\lambda-\lambda_0|\sim|\lambda |$ on the support of $a$; thus $B\les 1$ trivially.
This finishes the case $\xi >0>\xi '$.\\

\underline{Case 2:} To deal with the case $\xi >\xi '>0$, we use
(\ref{54}).  Thus,
$$f_-(\xi ',\lambda )=\alpha_-(\lambda )f_+(\xi ',\lambda )+\beta_-(\lambda )\overline{f_+(\xi ',\lambda )}$$
where
\begin{align*}
\alpha_-(\lambda )&=\frac{W(f_-(\cdot,\lambda ),\overline{f_+(\cdot,\lambda )})}{-2i\lambda}\\
\beta_-(\lambda )&=\frac{W(f_+(\cdot,\lambda ),f_-(\cdot,\lambda
))}{-2i\lambda}=\frac{W(\lambda )}{-2i\lambda}
\end{align*}
From our large $\lambda$ asymptotics of $W(\lambda )$ we deduce that
$\beta_-(\lambda )=1+O(\lambda ^{-1})$ and $\beta_-'(\lambda
)=O(\lambda^{-2})$. For $\alpha_-(\lambda )$ we calculate, again at
$\xi =0$,
\begin{align*}
W(f_-(\cdot,\lambda ),\overline{f_+(\cdot,\lambda
)})=\,&m_-(\xi,\lambda )(\overline{m}_+'(\xi,\lambda )
-2i\lambda\overline{m}_+(\xi,\lambda ))\\
&-\overline{m}_+(\xi,\lambda )(m'_-(\xi,\lambda )-2i\lambda m_-(\xi,\lambda ))\\
=\,&m_-(\xi,\lambda )\overline{m}'_+(\xi,\lambda)-m_-'(\xi,\lambda )\overline{m}_+(\xi,\lambda )\\
=\,&O(\lambda^{-1})
\end{align*}
so that
$$\alpha_-(\lambda )=O(\lambda^{-2})$$
with
$$\alpha'_-(\lambda )=O(\lambda^{-3}).$$
Thus, we are left with showing that
\begin{align}\label{70}
\sup_{\xi >\xi '>0}&\bigg|(\langle\xi\rangle\langle\xi '\rangle )^{-\frac{1}{2}}\int_{-\infty}^{\infty}e^{it\lambda^2}e^{i\lambda (\xi +\xi ')}\frac{\lambda (1-\chi (\lambda ))}{W(\lambda )}\alpha_-(\lambda )m_+(\xi,\lambda )m_+(\xi ',\lambda )\, d\lambda \bigg|
\les t^{-1}
\end{align}
and
\begin{align}\label{71}
\sup_{\xi >\xi '>0}&\bigg|(\langle\xi\rangle\langle\xi '\rangle )^{-\frac{1}{2}}\int_{-\infty}^{\infty}e^{it\lambda^2}e^{i\lambda (\xi -\xi ')}\frac{\lambda (1-\chi (\lambda ))}{W(\lambda )}\beta_-(\lambda )m_+(\xi,\lambda )\bar{m_+(\xi ',\lambda )}\, d\lambda \bigg|
\les t^{-1}
\end{align}
for any $t>0$.\\
This, however, follows by means of the exact same arguments which we
use to prove (\ref{64}).  Note that in (\ref{69}) the critical point
of the phase is
$$\lambda_0=-\frac{\xi +\xi '}{2t}$$
whereas in (\ref{70}) it is
$$\lambda_0=-\frac{\xi -\xi '}{2t}.$$
In either case it follows from $|\lambda_0|\gtr 1$ that $\xi\gtr t$.  Hence we can indeed argue as in case~1.\\
This finishes the proof of the lemma, and thus also Theorem
\ref{thm1}.
\end{proof}


\end{document}